\input amstex


\def\Aut{\text{Aut}}

\def\b1{\text{\bf 1}}

\def\CA{{\Cal A}}
\def\CB{{\Cal B}}
\def\CC{{\Cal C}}
\def\CD{{\Cal D}}

\def\CH{{\Cal H}}
\def\CI{{\Cal I}}

\def\CG{{\Cal G}}

\def\CO{{\Cal O}}

\def\CS{{\Cal S}}
\def\CT{{\Cal T}}
\def\CV{{\Cal V}}

\def\gr{\text{gr}}

\def\End{\text{End}}
\def\Hom{\text{Hom}}
\def\Sym{\text{Sym}}

\def\hotimes{\widehat{\otimes}}
\def\#{\,\check{}}

\def\fh{{\frak h}}

\def\fg{{\frak g}}

\def\fgl{{\frak{gl}}}
\def\fr{{\frak r}}
\def\fc{{\frak c}}

\def\id{\text{id}}
\def\tr{\text{tr}}

\def\pro{\text{pro}}

\def\Ker{\text{Ker}}

\def\Spec{\text{Spec}}

\def\holim{\text{holim}}
\def\dR{\text{dR}}

\def\\hat{\,\hat{}}

\def\pro{\text{pro-}}

\def\proab{\text{pro-}\!\CA b}
\def\limleft{\mathop{\vtop{\ialign{##\crcr
  \hfil\rm lim\hfil\crcr
  \noalign{\nointerlineskip}\leftarrowfill\crcr
  \noalign{\nointerlineskip}\crcr}}}}
\def\limright{\mathop{\vtop{\ialign{##\crcr
  \hfil\rm lim\hfil\crcr
  \noalign{\nointerlineskip}\rightarrowfill\crcr
  \noalign{\nointerlineskip}\crcr}}}}
\def\qlimleft{``\!\lim\!"}

\def\hra{\hookrightarrow}
\def\isol{\buildrel\sim\over\leftarrow}
\def\iso{\buildrel\sim\over\rightarrow} 

\def\lra{\longrightarrow}


\documentstyle{amsppt}
\NoBlackBoxes

\topmatter
\title     Relative continuous $K$-theory and cyclic homology    \endtitle
\author Alexander Beilinson  \endauthor
\leftheadtext{A.~Beilinson }
\address  Department of Mathematics, University of Chicago, Chicago, IL 60637, USA\endaddress
\email sasha\@math.uchicago.edu \endemail
\thanks The  author was supported in part by NSF	
Grant DMS-0401164.\endthanks
 \dedicatory To Peter Schneider
  \enddedicatory




\endtopmatter

\document


\head  Introduction. \endhead

\subhead{\rm 0.1}\endsubhead  
  Let  $A$ be a unital associative $p$-adic ring,  $p$ is a prime, so $A=\limleft A_i$ where $A_i :=A/p^i$.  Let $I\subset A$ be a two-sided ideal such that the $I$-adic topology equals the $p$-adic one, i.e.,  $p^m\in I $ and $  I^m \subset pA $ for large enough $m$.  Set $I_i :=IA_i$.
  
The projective system of rings $A_i$ and ideals $I_i$ yields the relative  $K$-theory  pro-spectrum  $K(A ,I)\\hat \, := \qlimleft   K(A_i ,I_i )$ and the cyclic chain pro-complex $CC(A)\\hat \, :=\qlimleft CC(A_i )$. 
Our main result, see 2.2, is a construction (subject to  minor  conditions on $A$) of a natural homotopy equivalence of the corresponding pro-spectra up to quasi-isogeny $$K(A,I)\\hat_{\Bbb Q} \iso CC(A)\\hat\, [1]_{\Bbb Q} . \tag 0.1.1$$

This homotopy equivalence is a continuous version of Goodwillie's isomorphism $K(R,J)\iso CC(R,J)[1]$  valid for any $\Bbb Q$-algebra $R$ and a nilpotent two-sided ideal $J$, see \cite{G} or \cite{Lo} 11.3. The construction follows closely that of Goodwillie, with a characteristic $p$ simplification (the difference between $ CC(A)\\hat $ and  $CC(A,I)\\hat \,$, being of torsion, disappears) and  the Malcev theory input replaced by a version of    Lazard's theory.

\subhead{\rm 0.2}\endsubhead Here is a geometric application (see 2.4):  Let $E$ be a $p$-adic field, $O_E$ its ring of integers,  $X$ be a proper  $O_E$-scheme with smooth generic fiber $X_E$,  $Y\subset X$ be any closed subscheme whose support equals the closed fiber. Set  $X_i := X\otimes\Bbb Z/p^i$. 
 Then (0.1.1) yields a canonical identification  of $\Bbb Q_p$-vector spaces  $$  \Bbb Q\otimes \limleft K^B_{n-1} (X_i , Y ) \iso \Bbb Q \otimes R\limleft  K^B_{n-1} (X_i ,Y)   \iso \oplus_a \, H^{2a -n}_{\dR}(X_E )/F^a . \tag 0.2.1$$ Here $K^B$ is the  non-connective $K$-theory (see \cite{TT} \S 6) and  $F^\cdot $ is  the Hodge filtration on the de Rham cohomology. 
 
The subjects left untouched: (i) comparison of the composition of (0.2.1) 
and the boundary map $K^B_{n}(Y)\to \limleft K^B_{n-1} (X_i , Y )$  with  the crystalline characteristic class map for $Y$  (here $Y$ can be arbitrary singular, see \cite{Cr} for the crystalline story in that setting), and (ii)  comparison with the $p$-adic regulators theory (see \cite{NN}).

\subhead{\rm 0.3}\endsubhead This note comes from an attempt to understand the work of Bloch, Esnault, and Kerz \cite{BEK} where an identification similar to (0.2.1) was constructed under extra assumptions on $X$.  
  The method of \cite{BEK} is  different (and it remains to be checked that the two isomorphisms coincide). Its input is McCarthy's identification of  K-theory of an {\it arbitrary} ring relative to a nilpotent ideal with the  relative
 {\it topological} cyclic homology, see \cite{DGM}. One has then  to identify the  topological cyclic homology  with  truncated de Rham cohomology. Unlike the classical cyclic homology case, this presents a  problem, and the passage taken in \cite{BEK} - with \'etale cohomology of $X_E$ as a way station - is not easy. 
 
It would be nice to find conditions on $(R,J)$, where
$p$ is nilpotent in $R$ and $J$ nilpotent, that would imply 
 $K_{i }(R,J)=H_{i-1} CC(R,J)$ for $i<p$, cf.~\cite{BEK} 8.5.

\subhead{\rm 0.4}\endsubhead I am  grateful to Spencer Bloch, H\'el\`ene Esnault, and Moritz Kerz whose ideas and interest were crucial for this work,  to  Mitya Boyarchenko, Volodya Drinfeld, John Francis, Lars Hesselholt, Dan Isaksen,  Dima Kaledin, Jacob Lurie,  Peter May, Chuck Weibel, and, especially, Nick Rozenblyum,  for the help, explanations, and/or exchange of letters, and to Linus Kramer for the tex assistance.

\head 1. Spectral preliminaries  \endhead

Given a sequence $\ldots \to A_1 \to A_0$ of abelian groups, then  $(\limleft A_i )\otimes \Bbb Q$  is not determined if we  know merely  $A_i \otimes \Bbb Q$, but we are in good shape if $A_i$ are known up to subgroups of torsion whose exponent is bounded by a constant independent of $i$. 
We will see  that some natural maps  known to be isomorphisms in rational homotopy theory are, in fact, invertible up to universally bounded denominators. Thus  they remain to be  isomorphisms in the rational homotopy theory of pro-spaces.

\subhead{\rm 1.1}\endsubhead {\it Spectra up to isogeny.} (a) Let $\CC$ be an additive category. We say that a nonzero integer $n$ {\it kills} an object $X$ of $\CC$ if $n\,\id_X =0$; $X$ is a {\it bounded torsion} object if it is killed by some  $n$ as above. A map $f :X\to Y$ in $\CC$ is an {\it isogeny} if there is  $g: Y\to X$    such that  $fg =n\,\id_Y$ and $gf =n\,\id_X$ for some nonzero integer $n$. 

We denote by $\CC \otimes \Bbb Q$ the category equipped with functor $\CC \to \CC\otimes\Bbb Q$, $X\mapsto X_{\Bbb Q}$, that is bijective on objects and yields identification 
 $\Hom (X_{\Bbb Q},Y_{\Bbb Q})=\Hom (X,Y)\otimes\Bbb Q $. We call $X_{\Bbb Q}$ the {\it object up to isogeny} that corresponds to $X$. Notice that $f$ is an isogeny if and only if
 $f_{\Bbb Q}$ is invertible, and $\CC \otimes \Bbb Q$ is the localization of $\CC$ with respect to isogenies; one has $X_{\Bbb Q}=0$, i.e., $X$ is isogenous to 0, if and only if $X$ is a bounded torsion object.

If $\CI$ is a category and $X, Y : \CI \to  \CC$ are functors, then a morphism of functors $f: X\to Y$ is said to be an isogeny if there is a morphism  $g: G\to F$    such that  $fg =n\,\id_Y$ and $gf =n\,\id_X$ for some nonzero integer $n$.  (If $\CI$ is essentially small, this means that $f$ is an isogeny in the category $\CC^{\CI}$.)

If $\CC$ is abelian then  $\CC \otimes \Bbb Q$ is abelian as well, and it coincides with the quotient  $\CC_{\Bbb Q}$  of $\CC$ modulo the Serre subcategory of bounded torsion objects. If $\CC$ carries a tensor structure then bounded torsion objects form an ideal, hence $\CC_{\Bbb Q}$ is a tensor category and $\CC \to \CC_{\Bbb Q}$ is a tensor functor.  In particular, the category
$\CA b$ of abelian groups yields the tensor abelian category  $\CA b_{\Bbb Q}$.

\remark{Remarks}  The functor $\CA b_{\Bbb Q}\to \CV ect_{\Bbb Q}$, $X_{\Bbb Q}\mapsto X\otimes \Bbb Q$, is {\it not} an equivalence of categories (for the category
 $\CV ect_{\Bbb Q}$  of $\Bbb Q$-vector spaces is the quotient of $\CA b$ modulo the Serre subcategory of groups all of whose elements are torsion). One has Ext$^1 (X_{\Bbb Q}, Y_{\Bbb Q})=$ Ext$^1 (X,Y)\otimes\Bbb Q$, 
 hence $\CA b_{\Bbb Q}$ has homological dimension 1 (since Ext$^1$ in $\CA b_{\Bbb Q}$ is right exact). $\CA b_{\Bbb Q}$ does  not admit infinite direct sums and products.
  \endremark

\enspace

It would be nice to refine the Quillen-Sullivan theory of rational homotopy types  to homotopy types up to isogeny (i.e., modulo Serre's class of abelian groups of bounded torsion). For the present article only stable theory is needed.

\enspace

(b) For an $\infty$-category treatment of the theory of spectra, see \cite{Lu2} 1.4.3. 

Let $\CS p$ be the stable $\infty$-category of spectra; by abuse of notation, we denote by $\CS p$ its homotopy category as well.
 This is a  t-category with heart $\CA b$: the subcategory $\CS p_{\ge 0}$ is formed by connective spectra,
the t-structure homology functor is the (stable) homotopy groups functor $X\mapsto \pi_\cdot (X)$, the t-structure truncations $X\to \tau_{\le n}X$ are  Postnikov's truncations. The t-structure is non-degenerate: if all $\pi_n (X)=0$ then $X=0$.
 $\CS p$ is a symmetric tensor $\infty$-category with respect to the smash-product $\wedge$, see \cite{Lu2}  6.3.2. It is right  t-exact and induces the usual tensor product on the heart $\CA b$; the unit object is the sphere spectrum $S$.

Let $\CD (\CA b)$ be the stable $\infty$-category of chain complexes of abelian groups, see  \cite{Lu2} 1.3.5.3, 8.1.2.8, 8.1.2.9;  we denote again by $\CD (\CA b)$ its homotopy category which is the derived category of $\CA b$ with its standard  t-structure. 
$\CD (\CA b)$ is a symmetric tensor $\infty$-category with the usual tensor product $\otimes$. 

The Eilenberg-MacLane functor 
EM: $\CD (\CA b )\to \CS p$ that identifies complexes with abelian spectra  is naturally  an $\infty$-category functor. EM is t-exact and equals   Id$_{\CA b}$ on the heart. It sends rings to rings and modules to modules, so
the  0th Eilenberg-Maclane spectrum $\Bbb Z_{\CS p}:=\text{EM} (\Bbb Z)$ is a unital ring spectrum and EM lifts to a functor $\CD (\CA b )\to \Bbb Z_{\CS p}$-mod (:= the $\infty$-category of $\Bbb Z_{\CS p}$-modules in $\CS p$). The latter functor is an equivalence of categories, see \cite{Lu2} 8.1.2.13.

Let  $\CS p^- :=\cup_n \,\CS p_{\ge 0}[-n]$ be the stable $\infty$-subcategory of eventually connective spectra; this is a tensor subcategory of $\CS p$. A map $X\to Y$ in $\CS p^-$ is called {\it quasi-isogeny} if all maps $\pi_n (X)\to \pi_n (Y)$ are isogenies, or, equivalently, all maps $\tau_{\le n}X \to \tau_{\le n}Y$ are isogenies in the homotopy category of spectra. Thus $X\in\CS p^-$ is  quasi-isogenous to 0 if each $\pi_n (X)$ is a bounded torsion group, or, equivalently,  every $\tau_{\le n} X$ is a bounded torsion spectrum; such $X$ form a thick subcategory. 

\proclaim{Lemma} This subcategory  is a $\wedge$-ideal in $\CS p^-$.
\endproclaim

\demo{Proof}  Suppose $X,Y\in \CS p^-$ and $Y$ quasi-isogenous to 0; let us check that $X\wedge Y$ is quasi-isogenous to 0. We can assume, replacing $X$, $Y$ by their shifts, that $X$ and $Y$ are connective. Then $\tau_{\le n}(X\wedge Y)=\tau_{\le n}(X\wedge \tau_{\le n}Y)$, and we are done. \qed
\enddemo

 Let $\CS p^-_{\Bbb Q}$ be the corresponding Verdier quotient category of $\CS p^-$; this is a symmetric tensor t-category. We call its objects {\it spectra up to quasi-isogeny}; for a spectrum $X$ we denote by $X_{\Bbb Q}$ the corresponding spectrum up to quasi-isogeny.  

Consider the derived categories $\CD^- (\CA b)$, $\CD^- (\CA b_{\Bbb Q})$ of bounded above complexes of abelian groups. Then  $\CD^- (\CA b_{\Bbb Q})$ is the  quotient of
$\CD^- (\CA b)$ modulo the thick subcategory of complexes with bounded torsion homology. EM sends $\CD^- (\CA b)$ to $\CS p^-$; passing to the  quotients, we get a t-exact functor 
$$\CD^- (\CA b_{\Bbb Q})\to \CS p^-_{\Bbb Q}. \tag 1.1.1$$

\proclaim{Proposition} This functor is an equivalence of tensor triangulated categories.
\endproclaim

\demo{Proof} Let $\Bbb Z_{\CS p}$-mod$^-$ be the category of eventually connective $\Bbb Z_{\CS p}$-modules and 
$\Bbb Z_{\CS p}$-mod$^-_{\Bbb Q}$ be its quotient  modulo objects which are quasi-isogenous to 0. The forgetful functor $\Bbb Z_{\CS p}$-mod$^- \to \CS p^-$ admits left adjoint  $X\mapsto   \Bbb Z_{\CS p}\wedge X$ which is a tensor functor. The adjoint functors $\Bbb Z_{\CS p}$-mod$^- \leftrightarrows \CS p^-$ yield, by passing to the quotients, the adjoint functors  $\Bbb Z_{\CS p}$-mod$_{\Bbb Q}^- \leftrightarrows \CS p_{\Bbb Q}^-$. By the above, it is enough to show that the latter functors are mutually inverse, i.e., that for $X\in \CS p^-$, $Y\in \Bbb Z_{\CS p}$-mod$^-$ the adjunction maps $a_X : X\to \Bbb Z_{\CS p}\wedge X$ and $a^\vee_Y :  \Bbb Z_{\CS p}\wedge Y \to Y$ are quasi-isogenies.

 One has $a_X = a_{S}\wedge \id_X$ where $a_S : S \to \Bbb Z_{\CS p}$ is the unit map. Since  $\pi_0 (a_S )=\id_{\Bbb Z}$, the groups $\pi_i (\CC one (a_S ))$ (the  stable homotopy groups of spheres)
are all finite, so $ \CC one (a_S )$, hence $\CC one (a_X )= \CC one (a_S)\wedge X $, is quasi-isogenous to 0, i.e., $a_X$ is a quasi-isogeny. Now $a^\vee_Y :  \Bbb Z_{\CS p}\wedge Y \to Y$  is a quasi-isogeny since $a^\vee_Y a_Y =\id_Y$. \qed
 \enddemo
\remark{Remark}  Let $c_n $ be an integer that kills $\tau_{\le n }\CC one (a_S)$ (e.g.~the product of exponents of the first $n$ stable homotopy groups of spheres). Then for every $X$ such that $\tau_{\le m}X=0$,   $c_{n-m}$ kills $\tau_{\le n}\CC one (a_X ) = \tau_{\le n}((\tau_{\le n-m}\CC one (a_S ))\wedge X)$. Therefore  $\pi_n (a_X ): \pi_n X \to \pi_n (\Bbb Z_{\CS p}\wedge X)$ is an isogeny of $\CA b$-valued functors on (every shift of) the category of connective spectra (see 1.1(a)).
\endremark
\enspace (c) We will need an $\CA b_{\Bbb Q}$-refinement of (a corollary of) the Milnor-Moore theorem:

\enspace

For a simplicial set (often referred to as topological space below) $P$ we denote by $C(P,\Bbb Z)$  the chain complex of $P$ and
 by $\bar{C}(P,\Bbb Z )$ the reduced chain complex which is
  the kernel of the augmentation map $C(P,\Bbb Z)\to \Bbb Z $, so for $p_0 \in P$ one has $\bar{C}(P,\Bbb Z )\iso C(P,\{ p_0 \};\Bbb Z)$ (:= the relative chain complex).
Suppose $P$ is connected. We say that $a\in H_n C (P,\Bbb Z)$  is {\it primitive} if $\Delta_* (a)\in H_n C (P\times P,\Bbb Z)=H_n (C(P,\Bbb Z)^{\otimes 2})$ equals $1\otimes a + a \otimes 1$, where 1 is the generator of $H_0 C (P,\Bbb Z)$ and $\Delta :P \hra P\times P$ is the diagonal map. 
We denote by Prim$ H_n C (P,\Bbb Z)$ the subgroup of primitive elements.  
 
 \remark{Exercises} (i) Show that the  projections $P\times P' \to P,P'$ yield 
isomorphisms  \newline Prim$ H_n C (P\times P',\Bbb Z)\iso\text{ Prim} H_n C (P,\Bbb Z)\oplus \text{Prim} H_n C (P',\Bbb Z)$. \newline
(ii) If $\pi_i (P)=0$ for $i>1$ and $\pi_1 (P)$ is abelian then $\text{ Prim} H_n C (P,\Bbb Z)=0$ for $n>1$.
\endremark

\enspace

 Denote by $\CT\! op_*$   the category of pointed topological spaces. We have the adjoint infinite suspension and infinite loop  functors $S^\infty : \CT\! op_* \leftrightarrows \CS p^- : \Omega^\infty$. For  $X\in \CS p^-$ and $P=(P, p_0 )\in \CT\! op_*$ let $b_X : S^\infty \Omega^\infty X \to X$ and $c_P : P \to \Omega^\infty S^\infty P $ be the adjunction maps, so one has  $\Omega^\infty (b_X) c_{\Omega^\infty X} =\id_{\Omega^\infty X}$ and $ b_{S^\infty P} S^\infty (c_P ) =\id_{S^\infty P}$. Let $h_P $ be the composition $P \to \Omega^\infty S^\infty P \to \Omega^\infty (\Bbb Z_{\CS p} \wedge S^\infty P) = \Omega^\infty \text{EM}(C(P,\{ p_0 \}; \Bbb Z ))$, the first arrow is $c_P$, the second arrow is $\Omega^\infty (a_{S^\infty P})$ (see the proof of the proposition in 1.1(b)). Then $\pi_n (h_P ): \pi_n (P,p_0 ) \to H_n C (P,\{p_0 \};\Bbb Z)=H_n C(P,\Bbb Z)$, $n\ge 1$, is the classical Hurewicz map; its image lies in Prim$ H_n C (P,\Bbb Z)$.

\proclaim{Theorem} For every connected spectrum $X$ the Hurewicz maps $\pi_n (h_{\Omega^\infty X} ) :\pi_n ( X)$ $=\pi_n (\Omega^\infty X, 0) \to \text{\rm Prim} H_n C (\Omega^\infty X,\Bbb Z)$  are isogenies. Moreover, their kernel and cokernel  are killed by  nonzero integers that depend only on $n$ (and not on $X$), i.e.,
$\pi_n (h_{\Omega^\infty} )$ are isogenies of $\CA b$-valued functors  on the category of connected spectra (see 1.1(a)). 
\endproclaim

\demo{Proof}\!\!\footnote{I am grateful to Nick Rozenblyum for the help with the proof.} 
(i)  The functor $X\mapsto \text{\rm Prim} H_n C (\Omega^\infty X,\Bbb Z)$ is additive: To see this, we need to check that for any maps of spectra $f, g : X\to Y$ and $\alpha\in \text{\rm Prim} H_n C (\Omega^\infty X,\Bbb Z)$ one has $(\Omega^\infty (f+g))_* (\alpha)=
(\Omega^\infty f)_* (\alpha )+ (\Omega^\infty g)_* (\alpha )$, which follows since
$\Omega^\infty (f+g)$ equals the composition $\Omega^\infty X \to \Omega^\infty X \times \Omega^\infty X \to  \Omega^\infty Y \times \Omega^\infty Y  \to \Omega^\infty Y$, the first arrow is $\Delta$, the second one is $\Omega^\infty f \times \Omega^\infty g$, the third one is the sum operation $+$.\footnote{ Use the fact that the restriction of $+$ to $\{0\}\times \Omega^\infty Y$ and $\Omega^\infty Y \times\{ 0\}$ is the identity map.}

\enspace

(ii) Consider a commutative diagram of spectra  $$\spreadmatrixlines{1\jot}
\matrix
S^\infty \Omega^\infty X
 &\to &\Bbb Z_{\CS p}\wedge
 S^\infty \Omega^\infty X   \\  
\downarrow  & & \downarrow  \\
X &\to &   \Bbb Z_{\CS p} \wedge X 
    \endmatrix \tag 1.1.2$$
 the horizontal arrows are $a_{S^\infty \Omega^\infty X}$ and $a_X$, the vertical arrows are $b_X$ and $ \id_{\Bbb Z_{\CS p}}\wedge b_X$.
It implies that  $\Omega^\infty (\id_{\Bbb Z_{\CS p}}\wedge b_X )h_{\Omega^\infty X}= \Omega^\infty (\id_{\Bbb Z_{\CS p}}\wedge b_X )\Omega^\infty (a_{S^\infty \Omega^\infty X}) c_{\Omega^\infty X} =\Omega^\infty (a_X) \Omega^\infty (b_X )     c_{\Omega^\infty X}=\Omega^\infty (a_X)$, hence $\pi_n (\id_{\Bbb Z_{\CS p}} \wedge b_X ) \pi_n (h_{\Omega^\infty X} ) = \pi_n (a_X )$. 
Let $$r_{nX} : \text{Prim}\, H_n C (\Omega^\infty X,\Bbb Z) \to \pi_n ( \Bbb Z_{\CS p} \wedge X ) \tag 1.1.3 $$ be the restriction of $\pi_n ( \id_{\Bbb Z_{\CS p}} \wedge b_X ): H_n 
C (\Omega^\infty X,\Bbb Z)\to \pi_n ( \Bbb Z_{\CS p} \wedge X ) $ to  primitive classes. Since $r_{nX}\pi_n (h_{\Omega^\infty X} )= \pi_n (a_X )$ is an isogeny of $\CA b$-valued functors on   connected spectra (see Remark in 1.1(b)), we see that $\pi_n (h_{\Omega^\infty })$ is an isogeny  if and only if  $ r_n$ is, and to check this it is enough to show that $\Ker\, r_{nX}$ is killed by a nonzero integer $e_n$ that does not depend on $X$.

\enspace

(iii) We first prove that $r_n$ is an isogeny on the subcategory of suspension spectra: Let us show  that $\Ker \, r_{n \,S^\infty P}$ for connected  $P=(P,p_0 )\in\CT\! op_*$  is killed by $n!$.
 
 As in \cite{M}, $\Omega^\infty S^\infty P$ identifies naturally with a free  $E_\infty$-space $F$ generated by $(P,p_0 )$. Let $\Gamma_\cdot$ be our $E_\infty$-operad. 
 Then $F$ is a union of closed subspaces $\{ 0\}=F_0 \subset F_1 \subset F_2 \subset \ldots $ where $(F_1, F_0 )=(P, \{ p_0 \})$  and $F_i$ is the image of the operad action map $ (\Gamma_i \times P^i )/\Sigma_i \to F$ where 
 $\Sigma_i$ is the symmetric group. Notice that $ (\Gamma_i \times P^i )/\Sigma_i $ is the homotopy quotient of $P^i$ modulo the action of  $\Sigma_i$.
 The pointed space $F_i /F_{i-1}$ equals $(\Gamma_i \wedge (P,p_0 )^{\wedge i})/\Sigma_i$ which is
 the homotopy quotient of $(P,p_0 )^{\wedge i} $ modulo  $\Sigma_i$, so $C(F_i ,F_{i-1};\Bbb Z)=
 C(\Sigma_i ,  C (P,\{ p_0 \}; \Bbb Z )^{\otimes i})$ where  $C(\Sigma_i ,\cdot )$ is the complex of group chains.  Thus the compositions   $C(\Sigma_i , \bar{C} (P,\Bbb Z )^{\otimes i}) \hra C(\Sigma_i , C (P,\Bbb Z )^{\otimes i}) \to C(F_i ,\{ 0\} ;\Bbb Z)$ split the filtration $C(F_i ,\{  0 \}; \Bbb Z )$ on $C(  F,\{  0 \};\Bbb Z)$ (see the beginning of 1.1(c) for the notation), so we have a canonical quasi-isomorphism   $$ \oplus_{i\ge 1}\,\, C(\Sigma_i , \bar{C} (P,\Bbb Z )^{\otimes i})\iso \bar{C}(F,\Bbb Z )= C (  F, \{ 0\} ,\Bbb Z ). \tag 1.1.4$$ 
 
The map  $\id_{\Bbb Z_{\CS p}}\wedge b_{S^\infty P} : \Bbb Z_{\CS p}\wedge S^\infty\Omega^\infty S^\infty P \to \Bbb Z_{\CS p}\wedge S^\infty P$  of abelian spectra amounts to a morphism $v_P : C(F,\{ 0\};\Bbb Z)\to C(P,\{ p_0\}; \Bbb Z)$ of complexes that can be described as follows.
 Consider
 the chain complex $C(P,\Bbb Z)$  as a simplicial abelian group freely generated by $P$, and $\bar{C}(P,\Bbb Z)$ as the kernel of the map of simplicial groups   $C(P,\Bbb Z)\to C(\{\text{point}\},\Bbb Z )$; ditto for $F$, etc. Then $v_P$ is  $\Bbb Z$-linear extension of the map of simplicial sets $F \to C(P,\Bbb Z)$ that extends the standard embedding $P\hra C(P,\Bbb Z)$ and transforms the $\Gamma_\cdot$-operations into addition.  Therefore $v_P$ is
the identity map on the first summand of (1.1.4) and it kills the rest of the sum. Thus  it remains to check that
 $H_n (\oplus_{i\ge 2}\, C(\Sigma_i , \bar{C} (P,\Bbb Z )^{\otimes i}))\cap\text{Prim} H_n (F ,\Bbb Z)$ is killed by $n!$. 
  
Set $\Phi^m := \oplus_{i\ge m}\, C(\Sigma_i , \bar{C} (P,\Bbb Z )^{\otimes i})$. We view $\Phi^\cdot$ as  a split filtration  on $C(F,\Bbb Z )$ via (1.1.4). 
Let us show that $m (H_n \Phi^m \cap \text{Prim} H_n (F ,\Bbb Z))\subset 
H_n \Phi^{m+1}$ for $m\ge 2$. Since $H_n \Phi^{n+1} =0$, this proves our assertion.

 We equip  $
  C(F\times F,\Bbb Z)=C(F,\Bbb Z )\otimes C(F,\Bbb Z )$ with the tensor product of the filtrations $\Phi^\cdot$. Since $\Delta_* : C(F,\Bbb Z )\to C(F\times F,\Bbb Z)$ is compatible with the operad product, $\Delta_*$  is compatible with the filtrations and $\gr^m_\Phi \Delta_* : C(\Sigma_m , \bar{C} (P,\Bbb Z )^{\otimes m})\to \oplus_{a+b=m} \,\, C(\Sigma_a , \bar{C} (P,\Bbb Z )^{\otimes a})\otimes C(\Sigma_{b} , \bar{C} (P,\Bbb Z )^{\otimes b})$ has $(a,b)$-component equal to the transfer map $\rho_{a,b}   : C(\Sigma_m , \bar{C} (P,\Bbb Z )^{\otimes m}) \to C(\Sigma_a \times \Sigma_{b}, \bar{C} (P,\Bbb Z )^{\otimes m})$ for the usual embedding $\Sigma_a \times\Sigma_b \hra \Sigma_m$. Now if $\alpha \in  \text{Prim} H_n (F ,\Bbb Z)$ lies in $  \Phi^m$, $m\ge 2$, then the image $\bar{\alpha} $ of $\alpha$ in $\gr^m_\Phi  H_n (F ,\Bbb Z)$ satisfies   $\rho_{a,b}(\bar{\alpha })=0$ if $a,b\ge 1$. In particular, $\rho_{1,m-1}(\bar{\alpha })=0$, and we are done since  $\Ker\,\rho_{1,m-1}$ is killed by $|\Sigma_m /\Sigma_{m-1}|=m$.

\enspace

(iv) To finish the proof, let us show that  $\Ker\, r_{nX}$ for any connected spectrum $X$ is killed by $e_n = n! c_n^{2n}$ where $c_n$ is as in Remark in 1.1(b).

To that end we will find  a map of spectra $s_{Xn} : \tau_{\le n} X\to \tau_{\le n} S^\infty \Omega^\infty X$ such that $\pi_i (\tau_{\le n}( b_X ) s_{Xn} )=\ell_n \id_{\pi_i (X)}$ for $i\le n$ where $\ell_n = c_n^2$. Assuming we have $s_{Xn}$, let us finish the proof.

One has $\tau_{\le n}( b_X ) s_{Xn} =\ell_n \id_{\tau_{\le n} X} + \epsilon$ where $\pi_i (\epsilon )=0$. Since the canonical filtration on $\tau_{\le n} X$ has length $n$, one has $\epsilon^n =0$. So, by (i),  $\tau_{\le n}( b_X ) s_{Xn}$ acts on Prim$H_n C(\Omega^\infty X ,\Bbb Z )$ as the sum of  the multiplication by $\ell_n$ map   and an operator $\epsilon$ such that $\epsilon^n =0$. The kernel $K$ of this action is killed by $\ell_n^{n}$ (for $\ell_n \id_K =-\epsilon|_K$). We are done since, by (iii) applied to $P= \Omega^\infty X$, the map $n! \tau_{\le n}(b_{X} ) s_{Xn}$ kills $\Ker\, r_{nX}=\Ker\, r_{n\, \tau_{\le n}X}$.\footnote{Notice that Prim$H_n C(\Omega^\infty X ,\Bbb Z ) \iso \text{Prim} H_n C(\Omega^\infty \tau_{\le n}X ,\Bbb Z )$.}

It remains to construct the promised $s_{Xn}$. For 
  any connected spectrum $Y$ let $q_Y :\tau_{\le n}(\Bbb Z_{\CS p} \wedge Y) \to \tau_{\le n} Y$ be a map of spectra  such that  $q_Y \tau_{\le n}(a_Y ) =c_n \id_{\tau_{\le n}Y} $.\footnote{Such a $q_Y$ exists since
 $c_n$ kills $\CC one (\tau_{\le n} ( a_Y ) )$.} 
  Let $u_i $, $i\in [1, n]$, be the composition $\pi_i (\Bbb Z_{\CS p}\wedge X )\to \pi_i (X)\to \pi_i (\Bbb Z_{\CS p}\wedge S^\infty \Omega^\infty X )$, the first arrow is $\pi_i (q_X )$, the second one is $\pi_i (h_{\Omega^\infty X})=\pi_i (a_{S^\infty \Omega^\infty X})\pi_i (c_{\Omega^\infty X} ) $. Let   $u:  \tau_{\le n} (\Bbb Z_{\CS p}\wedge X)\to \tau_{\le n} (\Bbb Z_{\CS p}\wedge S^\infty \Omega^\infty X)$ be a map  of spectra such that $\pi_i (u)=u_i$.\footnote{Our $u$ exists since $\tau_{\le n} (\Bbb Z_{\CS p}\wedge X)$ and $\tau_{\le n} (\Bbb Z_{\CS p}\wedge S^\infty \Omega^\infty X)$, being abelian spectra, are  homotopy equivalent to products of the Eilenberg-MacLane spectra corresponding to their homotopy groups.} 
   We define $s_{Xn}$ as the composition of maps $\tau_{\le n}X \to 
  \tau_{\le n} (\Bbb Z_{\CS p}\wedge X)$ $\buildrel{u}\over\to \tau_{\le n} (\Bbb Z_{\CS p}\wedge S^\infty \Omega^\infty X) \to  \tau_{\le n} (  S^\infty \Omega^\infty X)$ where the first arrow is $\tau_{\le n} (a_X)$ and the last arrow is $q_{S^\infty\Omega^\infty X}$. Finally, for $i\le n$ one has $\pi_i (\tau_{\le n}(b_X) s_{Xn} )  = \pi_i (b_X ) $ $  \pi_i (q_{S^\infty\Omega^\infty X} ) \pi_i (a_{S^\infty\Omega^\infty X} )\pi_i (c_{\Omega^\infty X} )\pi_i (q_X )\pi_i (a_X )=c_n^2 \pi_i (b_X)   \pi_i (c_{\Omega^\infty X} )= c_n^2 \id_{\pi_i (X)}  $, and we are done.
 \qed \enddemo

\remark{Remarks} (i) I do not know if the theorem remains true for connected H-spaces that are not infinite loop spaces.
\newline
(ii) The map $s_{Xn}$ from (iv) of the proof was constructed in an artificial manner;  it need not commute with morphisms of $X$. One can look for a natural $s_{Xn}$ (may be, with a better constant $\ell_n$). A possible approach: The ideal answer would be to find a natural section $s_X : X\to S^\infty \Omega^\infty X$ of $b_X$, and take  $s_{Xn}:= \tau_{\le n} (s_X )$ 
with $\ell_n =1$. Such an $s_X$ does not exist, but we have a natural section $c_{\Omega^\infty X}: \Omega^\infty  X\to \Omega^\infty S^\infty \Omega^\infty X$ of $\Omega^\infty (b_X )$ which   is {\it not}
an $E_\infty$-map for the standard $E_\infty$-structures on our spaces. Notice though that $\Omega^\infty S^\infty \Omega^\infty X$ is naturally an $E_\infty$-{\it ring} space with the product operation coming from the usual $E_\infty$-structure on $\Omega^\infty X$, and  $c_{\Omega^\infty X}$  is  
 an $E_\infty$-map for the product $E_\infty$-structure. A natural candidate for $s_X$ might be the logarithm of $c_{\Omega^\infty X}$. It is ill defined due   (at least) to
   non-integrality of the logarithm series, but  one can look for $n! \log \tau_{\le n}(c_{\Omega^\infty X})$ that would be an $E_\infty$-map for the standard $E_\infty$-structures, hence
 providing a natural $s_{Xn}$ with $\ell_n =n!$.
\endremark

\subhead{\rm 1.2}\endsubhead {\it The pro-setting.} 
(a) The basic reference for  pro-objects is   \cite{GV} \S 8. 

As in loc.~cit., one assigns to
 any category $\CC$  its category of pro-objects pro-$\CC$  together with a fully faithful functor  $\iota: \CC \hra $ pro-$\CC$, $X\mapsto   {}^\iota \! X$. The category pro-$\CC$ is closed under codirected  limits (i.e., for any $I^\circ$-diagram $Z_i$ in pro-$\CC$,   $I$  a   directed category, $\limleft Z_i$ exists), and $\iota$ is a universal functor from $\CC$ to such a category (i.e., every functor $\CC \to \CD$, where $\CD$ is closed under codirected limits, extends in a unique way to a functor pro-$\CC \to \CD$ that commutes with codirected limits). For $\CC $ small, the Yoneda embedding $\CC \to \text{Funct} (\CC  , \CS ets )^\circ$ extends naturally to a functor pro-$\CC \to \text{Funct} (\CC  , \CS ets )^\circ$; if $\CC$ is closed under finite limits, it identifies pro-$\CC$ with the subcategory of functors that commute with finite limits. Each pro-object 
 can be realized as $\qlimleft X_i :=\limleft {}^\iota \! X_i$ for
 some $I^\circ$-diagram $X_i$ in $\CC$, $I$ is directed, and
$\Hom (\qlimleft  X_i ,\qlimleft   Y_j )= \limleft_j \limright_i \Hom (X_i ,Y_j )$. If $\CC$ is abelian, then pro-$\CC$ is abelian and $\iota$ is exact. More generally, for any $I$ as above the functor $ \CC^I \to $ pro-$\CC$, $(X_i )\mapsto \qlimleft  X_i$, is exact. 
E.g.~we have abelian category pro-$\!\CA b$ of pro-abelian groups. The 
tensor product $\qlimleft  X_i \otimes \qlimleft   Y_j := \qlimleft  X_i \otimes Y_j $ makes it a symmetric tensor category.

\remark{Remarks} (i) If $\CC$ is  closed under arbitrary limits or colimits, then so is pro-$\CC$.
\newline (ii) For $\CC$  abelian, an object $\qlimleft X_i$ vanishes if and only if for every $i\in I$ one of the  maps $X_{i'} \to X_i$ of the diagram equals to zero.
\endremark

\enspace

(b) An $\infty$-category version of the story of \cite{GV} is explained in Chapter 5 of \cite{Lu1},  the $\infty$-category of pro-objects is treated in section 5.3.
 
Let pro-$\CS p$ be the $\infty$-category of pro-spectra, i.e., pro-objects of the $\infty$-category $\CS p$. This is a stable $\infty$-category by \cite{Lu2} 1.1.3.6. 
By abuse of notation, we denote its homotopy category by  pro-$\CS p$ as well. It carries a t-structure with the category of connective objects pro-$\CS p_{\ge 0}$ equal to 
pro-$(\CS p_{\ge 0})$. The heart of the t-structure equals
 $\proab$, the homology functor is $X=\qlimleft X_i \mapsto \pi_\cdot (X):= \qlimleft \pi_\cdot (X_i ) $, and the truncations are $\tau_{\ge n} \qlimleft X_i = \qlimleft \tau_{\ge n} X_i$, $\tau_{\le n} \qlimleft X_i = \qlimleft \tau_{\le n} X_i$. Let pro-$\CS p^- := \cup_n$pro-$\CS p_{\ge 0}[-n]$ be the t-subcategory  of eventually connective pro-spectra.  

The $\infty$-category space of morphisms of pro-spectra  $\CH om (\qlimleft X_i ,\qlimleft Y_j )$ is equal to holim$_j$hocolim$_i \CH om ( X_i ,Y_j )$;
passing to $\pi_0$ we get the set of morphisms in the homotopy category.

 By \cite{Lu2} 6.3.1.13, pro-$\CS p$ is a symmetric tensor $\infty$-category with tensor product $\qlimleft X_i \wedge \qlimleft Y_j= \qlimleft X_i \wedge Y_j $. The tensor product is right t-exact; it induces the usual tensor product on the heart pro-$\!\CA b$.

The above discussion can be repeated with ``spectrum" replaced by ``chain complex of abelian groups": we have stable $\infty$-category pro-$\CD (\CA b )$ of pro-complexes, its subcategory pro-$\CD (\CA b)^-$ of eventually connective  pro-complexes, etc. For a pro-complex $Z$ we denote by $H_\cdot Z$ its homology pro-abelian groups. We have the $t$-exact functor pro-EM:  pro-$\CD (\CA b )\to  $ pro-$\CS p$ which is identity on the hearts.

 A map $X\to Y$ in pro-$\CS p^-$ is 
 called {\it quasi-isomorphism}, resp.~{\it quasi-isogeny}, if the maps $\pi_n (X)\to \pi_n (Y)$ are isomorphisms, resp.~isogenies, for all $n$. So a pro-spectrum $X$ is quasi-isomorphic,  resp.~quasi-isogenous, to 0 if $\pi_n (X)$, resp.~ $\pi_n (X)_{\Bbb Q}$, vanish for all $n$. Such $X$ form thick subcategories of pro-$\CS p^-$ that are $\wedge$-ideals (repeat the proof of the lemma in 1.1(b)). We denote by pro-$\CS p^-_\pi$ and 
 pro-$\CS p^-_{\Bbb Q}$ the corresponding Verdier quotients; these are
 symmetric tensor t-categories with hearts $\proab$ and $\proab_{\Bbb Q}$,  the localizations pro-$\CS p^- \to$ pro-$\CS p^-_\pi \to $ pro-$\CS p^-_{\Bbb Q}$, $X\mapsto X_\pi \mapsto X_{\Bbb Q}$, are t-exact tensor functors. The discussion can be repeated for
pro-$\CD (\CA b)^-$.
 
 \remark{Exercise} For $X=\qlimleft X_i$ and $Y=\qlimleft Y_j$ in pro-$\CS p^-$ one has $\CH om\,  (X_\pi ,Y_\pi )=\holim_n   \CH om (\tau_{\le n}X,\tau_{\le n} Y)$, $\CH om (X_{\Bbb Q},Y_{\Bbb Q})=\holim_n \CH om (\tau_{\le n}X,\tau_{\le n} Y)\otimes \Bbb Q$.
\endremark 
 
 \remark{Remarks} (i) For any $Y\in \CS p^-$ 
 the pro-spectrum $X :=\qlimleft \tau_{\ge n}Y$ is quasi-isomorphic to zero;
 if $Y$ is not eventually coconnective, then
  $X\neq 0$.\footnote{I am grateful to Jacob Lurie for that example.} \newline
(ii) Probably pro-$\CS p^-_m$ coincides with
the eventually connective part of the  homotopy category of pro-spectra as defined in \cite{FI}.
\endremark

 \enspace
 
    Passing to the quotients, pro-EM yields the functor

$$\text{\rm pro-}\CD(\CA b)^-_{\Bbb Q}\to \text{\rm pro-}\CS p^-_{\Bbb Q}.\tag 1.2.1$$

 \proclaim{Proposition} This functor is an equivalence of tensor triangulated categories.\endproclaim
 
 \demo{Proof} Repeat the proof of the proposition in 1.1(b) using the remark after the proof to check that $\CC one (a_X )$ is quasi-isogenous to 0 for every    $X\in$ pro-$\CS p^-$.    \qed
 \enddemo

 (c) The theorem in 1.1(c) has an immediate pro-version as well. It will be used in \S  6 in the following manner:
 
 For   $X=\qlimleft X_i \in \pro\!\CS p^-$ one has a canonical map $\nu_X : C(\Omega^\infty X ,\Bbb Z)_{\Bbb Q} \to X_{\Bbb Q}$ defined as the composition $C(\Omega^\infty X ,\Bbb Z)= \Bbb Z_{\CS p}\wedge S^\infty \Omega^\infty X \to\Bbb Z_{\CS p}\wedge X  \leftarrow X$. Here $
 C(\Omega^\infty X ,\Bbb Z)=\qlimleft C (\Omega^\infty X_i ,\Bbb Z)$,
 the first arrow is $\id_{\Bbb Z_{\CS p}}\wedge b_X$, $b_X =\qlimleft b_{X_i}$, the second one is  quasi-isogeny $a_X =\qlimleft a_{X_i}$, see (1.1.2). 
 
Suppose $X$ is connected, i.e., $\pi_{\le 0}( X)=0$.  Let $V$ be a pro-complex such that $H_{\le 0}(V)=0$ and   $\psi : V_{\Bbb Q}\to C(\Omega^\infty X ,\Bbb Z)_{\Bbb Q}$ be a morphism such that for every $n$ the map  
 $H_n (\psi ): H_n V_{\Bbb Q}\to H_n C(\Omega^\infty X ,\Bbb Z)_{\Bbb Q}     $ is injective with image  Prim$H_n C(\Omega^\infty X ,\Bbb Z)_{\Bbb Q}$ $ := (\qlimleft $Prim$ H_n C(\Omega^\infty X_i ,\Bbb Z))_{\Bbb Q}$.

\proclaim{Proposition}  The map $\nu_X \psi  : V_{\Bbb Q} \to X_{\Bbb Q}$ is a quasi-isogeny.
 \endproclaim
 
 \demo{Proof} It is enough to check that $a_X  \nu_X \psi   : V \to \Bbb Z_{\CS p} \wedge X$ is a quasi-isogeny. One has  $H_n (a_X \nu_X \psi  )=r_{n X} H_n (\psi )$, see (1.1.3). We are done since $r_{nX}$ is an isogeny (see the proof of the theorem in 1.1(c)). \qed
 \enddemo
 
 (d) One has a natural functor of stable $\infty$-categories holim: pro-$\CS p  \to \CS p$, $\qlimleft X_i \mapsto \holim X_i$,  which is left t-exact. Its restriction to the  subcategory pro$^{\aleph_0}$-$\CS p^-$ of those $ \qlimleft X_i$ that the set  $I$ of indices $i$ is countable, has homological dimension 1,\footnote{For a precise homological dimension assertion if $I$ is arbitrary, see \cite{Mi}.} so we have triangulated functor holim:  pro$^{\aleph_0} $-$\CS p^-  \to \CS p^-$ and, passing to the quotient categories, holim:  pro$^{\aleph_0}$-$\CS p^-_{ \Bbb Q} \to \CS p^-_{\Bbb Q}$. A similar functor  holim: $\text{pro}^{\aleph_0}$-$\CD (\CA b )^- \to \CD(\CA b)$  is the right derived functor of the projective limit functor $\qlimleft C_i \mapsto \limleft C_i$, and the two holim's are compatible with the (pro-) EM functors.
   
 \head 2. The main theorem and a geometric application  \endhead

\subhead{\rm 2.1}\endsubhead For  an associative unital $V$-algebra $R$, $V$ is a commutative ring, we denote by $C(R/V)$ the relative
 Hochschild complex (see \cite{Lo} 1.1.3), so $C (R/V)_n =R^{\otimes n+1}$ for $n\ge 0$ (here $\otimes =\otimes_V$), $C (R/V)_n = 0$ otherwise. If $R$ is $V$-flat, then $C(R/V)$ equals $R\otimes^L_{R\otimes_V \! R^\circ} R$. Let $CC(R/V)$ be the cyclic complex as in \cite{Lo} 2.1.2: this is the
 total complex of the cyclic bicomplex $CC(R/R)_{\cdot ,\cdot}$; here
$CC(R/V)_{m,n}=R^{\otimes n+1}$ for $m,n\ge 0$ and $C(R/V)_{m,n}=0 $
otherwise; the $n$th row complex  is the chain complex for the $\Bbb Z/n+1$-action $r_0 \otimes \ldots\otimes r_n \mapsto (-1)^n r_n \otimes r_0 \otimes \ldots \otimes r_{n-1}$
on $R^{\otimes n+1}$; the even column complexes equal $C(R/V)$ and the odd ones are acyclic. We denote by $\Phi_m$ the increasing filtration on $CC(R/V)$ that comes from the filtration $CC(R/V)_{\le 2m,\cdot}$ on the bicomplex; thus $\Phi_{-1}=0$ and the projection $\gr_m^\Phi CC(R)\twoheadrightarrow CC(R/V)_{2m,\cdot}= C(R/V)[2m]$, $m\ge 0$, is a quasi-isomorphism.

Apart from 2.3--2.4, we consider only Hochschild and cyclic complexes with $V=\Bbb Z$ and denote them simply $C(R)$, $CC(R)$.

\subhead{\rm 2.2}\endsubhead {\it The main theorem.} Let $A$, $I$  be as in 0.1;   as in loc.~cit., $A_i := A/p^i$ and $I_i $ is the image of $I$ in $A_i$. Consider the relative $K$-theory
pro-spectrum $ K (A,I)\\hat  \, := \qlimleft  K(A_i ,I_i )$ and the cyclic homology pro-complex $CC (A)\\hat \, :=\qlimleft  CC (A_i )$ (see 3.1 below for a reminder and the notation).

\proclaim{Theorem} Suppose $A$ has bounded $p$-torsion\footnote{I.e., the $p$-torsion subgroup of $A$ is killed by some $p^n$.} and $A_1$ has finite stable rank (see \cite{Su}). 
Then there is a natural quasi-isogeny  $$ K(A,I)\\hat_{\Bbb Q}\iso CC (A)\\hat\, [1]_{\Bbb Q}  . \tag 2.2.1$$
\endproclaim

The proof takes \S\S 3--6. The outline of the construction, which is a continuous version of Goodwillie's story \cite{G} or \cite{Lo} 11.3.2, is as follows:

- Consider the pro-complex $C(\fgl_r (A)\\hat \,  ):=\qlimleft C(\fgl_r (A_i ))$  of Lie algebra chains; here $\fgl_r$ is the
Lie algebra of matrices.
By a version of the Loday-Quillen-Tsygan theorem \cite{Lo} 10.2.4, there is a canonical map of pro-complexes $C(\fgl_r (A)\\hat \, )\to \Sym (CC(A)\\hat \,[1])$ which is a quasi-isogeny in degrees $\le r$.

- Let $GL_r (A)^{(m)}:=\Ker (GL_r (A)\to GL_r (A_m )) $ be the congruence subgroup of  level  $m\ge 2$, and
 $C (GL_r (A) ^{(m)} \\hat  \, ) := \qlimleft C( GL_r (A) ^{(m)} /  GL_r (A) ^{(m+i)} ,\Bbb Z )$ be the pro-complex of its
 group chains. A version of Lazard's theory from Ch.~V of \cite{L} provides  (we use the bounded $p$-torsion condition here)  a
canonical quasi-isogeny  $C (GL_r (A)^{(m)}   \\hat  \, )_{\Bbb Q}\iso C(\fgl_r (A)\\hat \, )_{\Bbb Q}$.

- We realize $ K(A,I)\\hat   \, $ using Volodin's construction.  By Suslin's results \cite{Su} (we use the finiteness of stable rank condition here) one can compute the relative $K_i$ pro-groups using $GL_r$'s with bounded $r$ depending on $i$. Viewed up to isogeny, they coincide with the primitive part of  $H_i C(GL_r (A)^{(m)} \\hat\, )  $.

Combining the above quasi-isogenies and   using 1.2(c), we get (2.2.1).

\subhead{\rm 2.3}\endsubhead In the next subsection we explain a geometric application of the main theorem. We need three technical lemmas; the reader can skip them at the moment to return when necessary.

(a) Let $R$ be a ring and $I$ its two-sided ideal. Let $K^B$ be the non-connective relative $K$-theory spectrum from  \cite{TT} \S 6, so $K^B_n (R) =K_n (R)$ for $n\ge 0$ and there is a long exact sequence $\ldots \to K^B_n (R,I)\to K^B_n (R)\to K^B_n (R/I)\to K^B_{n-1}(R,I)\to \ldots$

\proclaim{Lemma} If $I$ is nilpotent then $K^B (R,I)$ is connected, i.e., $K^B_n (R,I)=0$ for $n\le 0$.
\endproclaim

\demo{Proof} The assertion of the lemma amounts to surjectivity of the map $K_1 (R)\to K_1 (R/I)$ and bijectivity of the maps $K^B_n (R)\to K^B_n (R/I)$ for all $n\le 0$. 

 Every idempotent in Mat$_n (R/I)$ lifts to an idempotent in Mat$_n (R)$, and the map Isom$(P,P')\to $ Isom$(P/IP, P'/IP')$ is surjective for all projective $R$-modules $P$, $P'$. Therefore $K_1 (R)\twoheadrightarrow K_1 (R/I)$ and $K_0 (R)\iso K_0 (R/I)$.

Induction by $-n$: Suppose $K^B_n (R) \iso K^B_n (R/I)$ for all $(R,I)$ with $I$ nilpotent. Then $K^B_{n-1} (R) \iso K^B_{n-1} (R/I)$. Indeed, the map from the canonical exact sequence $0\to K^B_n (R)\to K^B_n (R[t]) \oplus K^B_n (R[t^{-1}])\to K^B_n (R[t,t^{-1}])\to K^B_{n-1}(R)\to 0$ of \cite{TT} 6.6 to the similar exact sequence for $R/I$  is an isomorphism at each term but the last one by the induction assumption. So it is an isomorphism everywhere, q.e.d.     \qed
\enddemo

(b) Let $X\\hat$ be a formal scheme, so $X\\hat\,$ is direct limit of a directed family of closed subschemes $X_i$ where the ideals of the embeddings $X_i \hra X_j$ are nilpotent. Let $Y$ be a closed subscheme of  $X\\hat$ such that
 $X_i$ contain $Y$ and the ideal of $Y$ in $X_i$ is nilpotent. One has the projective system of the non-connective relative $K$-theory spectra $K^B (X_i ,Y)$. 

\proclaim{Lemma} If $Y$ is quasi-compact and separated then 
the pro-spectrum $\qlimleft K^B (X_i ,Y)$ is eventually connective. \endproclaim

\demo{Proof} Pick a finite open affine Zariski covering $\{ Y_\alpha \}_{\alpha \in \CA}$ of $Y$; let $\{ X_{i\alpha    }\}$ be the corresponding covering of $X_i$. Let $\frak S$ be the set of non-empty subsets $S$ of $\CA$ ordered by inclusion. We have an $\frak S$-diagram $S\mapsto Y_S :=\cap_{\alpha\in S} Y_\alpha$
 of affine schemes and open embeddings, and similar $\frak S$-diagrams $S \mapsto X_{iS}$. By  \cite{TT} 8.4, one has  $K^B (X_i ,Y)\iso  \holim_\frak S \, K^B (X_{iS } ,Y_{S})$. Thus  $K^B (X_i ,Y)$ has a finite filtration with $\gr^n K^B (X_i ,Y)=\Pi_{|S|=n} K^B (X_{iS } ,Y_{S})[1-n]$.
 Since $X_{iS}$ is affine and the ideal of $Y_S$ in $X_{i S}$ is nilpotent, the above lemma shows that the spectrum  $K^B (X_{iS } ,Y_{S})$ is   connected.  Therefore $\pi_{-n} K^B (X_i ,Y)=0$ for $n\ge |\CA|-1$, and we are done.  \qed
\enddemo

(c) Let  $k$ be a perfect field of characteristic $p$, $W=W(k)$ be the Witt vectors ring, and $R$ be a unital associative flat $W_i $-algebra where $W_i :=W/p^i$. 

\proclaim{Lemma} The evident map  $\tau  : CC(R)\twoheadrightarrow CC(R/W_i )$  is a quasi-isomorphism.
\endproclaim
\demo{Proof} The map $\tau $ is compatible with the filtrations $\Phi$ (see 2.1), so it is enough to check that $\tau $ is a filtered quasi-isomorphism, i.e.,  that the map of   Hochschild complexes $ C(R) \to C(R/W_i )$ is a quasi-isomorphism.  Both   $C(R)$ and $C(R/W_i )$ are $\Bbb Z/p^i$-flat (since $R$ is $W_i$-flat), so it is enough to check that the map $C(R)\otimes\Bbb Z/p  \to C(R/W_i )\otimes\Bbb Z/p$ is a quasi-isomorphism. The latter  is the map $C(R_1 )\to C(R_1 /k )$, $R_1 :=R/pR$, so, replacing $R$ by $R_1$, we can assume that $i=1$. 

We check first that $C(k)=C(k/\Bbb F_p ) \to   k$ is a quasi-isomorphism. It is enough to consider the case when $k$ is the perfectization of a field $k'$ finitely generated over $\Bbb F_p$. Then $k'$ is a separable extension of a purely transcendental extension of $\Bbb F_p$, and $H_i C(k' /\Bbb F_p ) \iso \Omega^i (k'/\Bbb F_p )$ by \cite{Lo} 3.4.4. Since Frobenius kills
$\Omega^{>0} (k'/\Bbb F_p )$, one has $H_{>0} C(k)=0$, q.e.d.

If $R$ is an arbitrary $k$-algebra, then let $P$ be an $R\otimes_{\Bbb F_p}R^\circ$-flat resolution of $R$. The terms of $P$ are $k\otimes_{\Bbb F_p}k$-flat, hence $P_k := P \otimes_{k\otimes_{\Bbb F_p}k} k$ is again a resolution of $R$ (indeed, $P_k \iso R\otimes^L_{k\otimes k}k= R\otimes_k (k\otimes^L_{k\otimes k}k)\iso R\otimes_k C(k)$ which is $R$ by the above). Thus $P_k$ is an $R\otimes_k R^\circ$-flat resolution of $R$, hence $C(R)\iso P\otimes_{R\otimes_{\Bbb F_p}R^\circ}R = P_k \otimes_{R\otimes_{k}R^\circ}R \iso C(R/k)$, and we are done.
\qed
\enddemo

\subhead{\rm 2.4}\endsubhead Let $E$ be a $p$-adic field, $O_E$ be its ring of integers, so $O_E$ is a complete mixed characteristic dvr with perfect residue field $k$ of characteristic $p$. We denote by  $O_E\text{-mod}$  the category of $O_E$-modules, by $\CD^- (O_E\text{-mod})$ the derived category of bounded above complexes, and by $\CD_{\text{fg}}^- (O_E\text{-mod})$ its full subcategory formed by  complexes whose homology are finitely generated $O_E$-modules.
Consider the functor $\CD^- (O_E\text{-mod})\to $ pro-$\CD  (\CA b)^-$, $M\mapsto M\hotimes\Bbb Z_p  :=  \qlimleft M\otimes^L \Bbb Z/p^i =\qlimleft \CC one (M\buildrel{p^i}\over\to M)  $. Notice that for $M\in \CD_{\text{fg}}^- (O_E\text{-mod})$
 the evident map $M\to \holim_i \, M\otimes^L \Bbb Z/p^i$ is a quasi-isomorphism  (which is enough to check  for $M=O_E$), and the restriction of our functor $\cdot \hotimes\Bbb Z_p$ to $\CD_{\text{fg}}^- (O_E\text{-mod})$ is t-exact 
(since for any  finitely generated $O_E$-module $M$ one has $\qlimleft \text{Tor}_1 (M,\Bbb Z/p^i )=0$).

Let $\CD_{\text{fd}}^- (\text{Vect}_E )$ be the derived category of bounded above complexes of $E$-vector spaces with finite-dimensional homology. The functor  $\CD_{\text{fg}}^- (O_E\text{-mod})\to \CD_{\text{fd}}^- (\text{Vect}_E )$, $M\mapsto M\otimes\Bbb Q$, identifies the target with the Verdier quotient of $\CD_{\text{fg}}^- (O_E\text{-mod})$ modulo   thick subcategory of complexes with torsion cohomology. The functor $\cdot\, \hotimes \Bbb Z_p$ sends the latter subcategory to pro-complexes quasi-isogenous to zero, hence it yields a t-exact functor $$\CD_{\text{fd}}^- (\text{Vect}_E)\to \text{\rm pro-}\CD (\CA b)^-_{\Bbb Q} \tag 2.4.1 $$ which is evidently faithful (fully faithful if $E=\Bbb Q_p$). 

 Let $X$ be a proper  $O_E$-scheme with smooth generic fiber $X_E$.  
  Let $R\Gamma_{\!\text{dR}}(X_E ):=R\Gamma (X_E, \Omega^\cdot_{X_E /E} )$ be the de Rham  complex of $X_E$ and $F^\cdot$ be its Hodge filtration. The complexes $R\Gamma_{\!\text{dR}}(X_E )/F^a$ are bounded and have finite-dimensional cohomology.

Let $Y\subset X$ be a closed subscheme whose support equals the closed fiber. Set  $X_i := X\otimes\Bbb Z/p^i$, so $Y\subset X_i$ for large enough $i$. By   2.3(b) we have the   pro-spectrum  $K^B (X,Y)\\hat\, := \qlimleft K^B (X_i ,Y)\in $  pro-$\CS p^-$.

 \proclaim{Theorem} There is a natural quasi-isogeny of pro-spectra $$K^B (X,Y)\\hat_{\Bbb Q} \iso \oplus_a \,  (R\Gamma_{\!\text{dR}}(X_E )/F^a ) [2a-1] . \tag 2.4.2$$ Here  the target is seen  via (2.4.1) as an object  of {\rm pro-}$\CD  ( \CA b)^-_{\Bbb Q}$, hence a pro-spectrum up to quasi-isogeny.
 \endproclaim
 
Applying $\holim$ from 1.2(d), we get identification (0.2.1).

\demo{Proof} (i) Pick  a finite open affine covering $\{X_\alpha \}$
of $X$; let $\{ Y_\alpha := Y\cap X_\alpha \}$ be the induced covering of $Y$. As in the proof of the  lemma in 2.3(b), we write $X_S :=\cap_{\alpha \in S} X_\alpha =\Spec\, R_S$, $Y_S =\Spec\, R_S /I_S$  for $S\in \frak S$, etc. 
By loc.~cit.~and the lemma in 2.3(a),   $K^B (X,Y)\\hat =\qlimleft \holim_{\frak S} K(X_{i S}, Y_S )=   \holim_{\frak S} \, K(X_{ S}, Y_S )\\hat\,$. Passing to pro-spectra up to isogeny, we get $K^B (X,Y)\\hat_{\Bbb Q} = \holim_{\frak S} \, K(X_{S}, Y_{S})\\hat_{\Bbb Q}$. 

Set $CC (X_{ i}):=  \holim_{\frak S} \,CC (X_{iS })$, 
$ CC(X_S )\\hat := \qlimleft   CC (X_{i S })$, and
$CC(X)\\hat := \qlimleft  CC (X_{ i})= \holim_{\frak S} \,  CC(X_S )\\hat\,$.    The rings
$A_S =\limleft R_{iS}$ satisfy the conditions of the
theorem in 2.2, so (2.2.1) applied to 
 $A_S $ and   ideals $I_S A_S$  provides then a canonical identification 
 $$K^B (X,Y)\\hat_{\Bbb Q} \iso   CC(X )\\hat_{\Bbb Q}[1]. \tag 2.4.3 $$  
  
(ii)  The terms of (2.4.2) do not change if we replace $E$ by the field of fractions of $W=W(k)$, so we can assume that $O_E =W$. 

Let $X_{\text{red}} \subset X$ be the reduced scheme. Since $X_E$ is smooth and $X$ has finite type, $\CO_{X_{\text{red}}}$ is the quotient of $\CO_X$ modulo ideal killed by  high enough power of $p$, so the maps
$CC(X_{iS })\\hat \to CC (X_{\text{red}\, iS })\\hat $ are quasi-isogenies. Therefore, by (2.4.3),  $K^B (X,Y)\\hat_{\Bbb Q} \iso K^B (X_{\text{red}},Y_{\text{red}})\\hat_{\Bbb Q}$. Thus the terms of (2.4.2) do not change if we replace $X$ by $X_{\text{red}}$, so we can assume that $X$ is flat over $W $.

\enspace

(iii)  Consider  $W$-complexes
 $CC (X_S /W):= CC(R_S /W)$; they are $W$-flat since $R_S$ is $W$-flat. 
 Set $CC(X/W):= \holim_{\frak S}\, CC (X_S /W)$. 
 By  2.3(c) one has  $CC (X_{iS })\iso CC(X_{iS }/W_i )=CC (X_S /W)\otimes \Bbb Z/p^i $. Therefore $$CC(X  )\\hat \, \iso \qlimleft  CC(X /W)\otimes \Bbb Z/p^i . \tag 2.4.4$$
 
 The complex $CC (X/W )$ lies in $\CD_{\text{fg}}^- (W\text{-mod})$, i.e., the homology of $ CC (X/W )$ are finitely generated $W$-modules: Indeed, since 
filtration $\Phi $ on $CC$ (see 2.1) is finite on every homology group and $\CS$ is finite,
 it suffices to show that $\gr^{\Phi}_m  CC (X/W )$ $:=  \holim_{\frak S}\, \gr^{\Phi}_m CC(X_S /W)\in   \CD_{\text{fg}}^- (W\text{-mod})$. Since 
$ \gr^{\Phi}_m CC(X_S /W)\iso   C(X_S /W )[2m]$, it is enough
to check that
$ \holim_{\frak S}\, H_n C(X_S /W )\in  \CD_{\text{fg}}^- (W\text{-mod}) $. Now $ H_n C(X_S /W )$  $=\Gamma (X_S, \CH_n )$, $\CH_n := H_n  L\Delta^* \Delta_* \CO_X $ where  $\Delta : X\hra X\times_{\Spec\, W}X$ is the diagonal embedding. Since $\CH_n$ is a coherent $\CO_X$-module and 
 $X/W$ is proper, the homology of $ \holim_{\frak S}\, H_n C(X_S /W )=R\Gamma (X,\CH_n )$ are  finitely generated $W$-modules, q.e.d.
 
 \enspace

(iv) Combining (2.4.3) and (2.4.4) we get a canonical identification $$K^B (X,Y)\\hat_{\Bbb Q}\iso (\qlimleft  CC(X /W)\otimes \Bbb Z/p^i )_{\Bbb Q}[1]. \tag 2.4.5$$ By (iii) and the definition of (2.4.1),
to get (2.4.2) it remains to produce  a natural quasi-isomorphism $$CC(X,W)\otimes\Bbb Q \iso  \oplus_a \,  (R\Gamma_{\!\text{dR}}(X_E )/F^a ) [2a-2]. \tag 2.4.6$$ Now $CC(X/W)\otimes\Bbb Q =\holim_{\frak S}\, CC(X_{S E}/E )$, and (2.4.6) comes from a canonical  quasi-isomorphism $CC(R/E ) \iso \oplus_a \, (\Omega^\cdot (R/E ) /F^a )[2a-2]$ from \cite{Lo}  3.4.12\footnote{By loc.cit., this quasi-isomorphism is defined as the
composition of quasi-isomorphisms \newline $CC(R/E ) \hookleftarrow \CB (R/E ) \twoheadrightarrow \oplus_a \, (\Omega^\cdot (R/E ) /F^a )[2a-2]$
 from \cite {Lo} 2.1.7, 2.1.8 and \cite{Lo} 2.3.6, 2.3.7.} valid for any smooth commutative algebra $R$ over $E\supset \Bbb Q$, and applied to
algebras $R=R_S \otimes \Bbb Q$, and we are done.
 \qed \enddemo

\head 3. The Loday-Quillen-Tsygan isomorphism  \endhead

The
Loday-Quillen-Tsygan theorem, see \cite{Lo} 10.2.4,  provides a canonical quasi-isomorphism $C (\fgl (R))\iso \Sym (CC (R)[1])$ for any $\Bbb Q$-algebra $R$. We adapt the  argument of loc.~cit.~for our continuous setting.

\subhead{\rm 3.1}\endsubhead For an associative unital ring $R$ let 
$ C^\lambda (R)$ be the Connes complex (see
 \cite{Lo} 2.1.4). This  is  the quotient complex of the Hochschild complex $C  (R)$: namely, $C^\lambda (R)_n$ is the coinvariants of the $\Bbb Z/n+1$-action on $R^{\otimes n+1}$ (see 2.1).
There is an evident projection $\pi_R : CC  (R)\twoheadrightarrow C^\lambda  (R)$  that equals the projection $C  (R)
\twoheadrightarrow C^\lambda  (R)$ on $CC(R)_{0,\cdot}$ and kills the rest of the bicomplex $CC(R)_{\cdot ,\cdot}$.

\proclaim{Lemma}  The homology group $H_n \Ker\, (\pi_R )$ is killed by $n!$.
\endproclaim

\demo{Proof} Consider the increasing filtration by the bicomplex row number on $CC (R)_\cdot$ and its subcomplex $\Ker (\pi_R )$. Then $H_n \gr_m \Ker (\pi_R )$ equals $H_{n-m}(\Bbb Z/m+1 , R^{\otimes m+1})$ for $n>m>0$ and is 0 otherwise. It is killed by   $m+1$, hence the assertion. \qed
\enddemo

\subhead{\rm 3.2}\endsubhead Let $A$ be any $p$-adic unital ring, so $A=\limleft A_i$ where $A_i :=A/p^i$. 
We have pro-complexes $C  (A)\\hat \, :=\qlimleft  C(A_i )$ and, similarly,  
 $C^\lambda    (A)\\hat \,    $, $CC  (A)\\hat \, $. 
 Consider the projection $\pi_{A\\hat }: CC (A)\\hat \, 
\twoheadrightarrow  C^\lambda    (A)\\hat \, $.

\proclaim{Corollary}  $ \pi_{A\\hat} $ is a quasi-isogeny.
\endproclaim

\demo{Proof} Consider the pro-complex $\Ker (\pi_{A\\hat})$. The pro-group $H_n \Ker (\pi_{A\\hat})$ is equal to $\qlimleft  H_n  \Ker (\pi_{A_i})$, so it is killed by
 $n!$   according to the above lemma.
We are done by the exact triangle $ \Ker (\pi_{A\\hat } )\to CC (A)\\hat \, \to C^\lambda   (A)\\hat \, $. \qed
\enddemo

 \subhead{\rm 3.3}\endsubhead For a Lie algebra $\fg$ we denote by $C  (\fg )$ its  Chevalley complex with trivial coefficients (see \cite{Lo} 10.1.3), so  $C(\fg )_n =\Lambda^n \fg $ (the exterior power over $\Bbb Z$) and
$H_n (\fg):= H_n C  (\fg )$ are the Lie algebra homology groups. Our 
$C  (\fg )$ is a cocommutative counital dg coalgebra, the coproduct map $\delta: C  (\fg )\to C  (\fg \times  \fg)=C(\fg)\otimes C  (\fg)$  comes from the diagonal map $\fg\to \fg \times \fg $.

 Consider the symmetric dg algebra $\Sym ( C^\lambda (R)[1])$. This is a commutative and cocommutative (unital and counital) Hopf dg algebra; the coproduct comes from the
 diagonal map  $C^\lambda (R)[1]\to C^\lambda  (R)[1]\times  C^\lambda  (R)[1]$. 
 
We recall the  key construction of Loday-Quillen and Tsygan, cf.~\cite{Lo} 10.2. Consider the matrix Lie algebra $\fgl_r (R)$, $r\ge 1$.

 \proclaim{Theorem-construction} There is  a canonical  morphism of dg coalgebras
 $$\kappa_R =\kappa^r_{R}  :   C  (\fgl_r (R)) \to \text{ \rm{Sym}} ( C^\lambda (R)[1]). \tag 3.3.1$$ The maps $C (\fgl_r (R))_n \to  (\text{\rm Sym} ( C^\lambda (R)[1]))_n$ are surjective for $n\le r $. 
\endproclaim

\demo{Proof} 
We need a preliminary construction from \cite{Lo} 9.2. Let $V$ be a free abelian group of rank $r$ and $I$ be a finite set of order $n$.
The symmetric group $\Sigma_I$ acts on $V^{\otimes I}$ transposing the factors. Let $\theta^*_I :  \Bbb Z [\Sigma_I ] \to  \End (V)^{\otimes I}$ be the composition
 $ \Bbb Z [\Sigma_I ] \to \End (V^{\otimes I})\buildrel{\sim}\over\leftarrow \End (V)^{\otimes I}$ of the $\Sigma_I$-action map and the inverse to the tensor product map. The source and   target of $\theta^*_I$ are naturally self-dual: the elements $\sigma$ of $\Sigma_I$ form an orthonormal base of the source,  
 the duality for the target is $f,g \mapsto \tr (fg)$. Thus we have the dual map $\theta_I =\Sigma\,  \theta_\sigma \sigma:  \End (V)^{\otimes I}\to  \Bbb Z [\Sigma_I ]$. It commutes with  the action of $\Aut (V) \times \Sigma_I$; here  $\Sigma_I$ acts on the target by conjugation.
  
  An explicit formula for $\theta_\sigma$: Let $\{ I_\alpha \}$ be the orbits of $\sigma$. So $I=\sqcup I_\alpha$ and $I_\alpha$ are naturally cyclically ordered: we  write $I_\alpha = \{ i^\alpha_1 ,\ldots ,i^\alpha_{n_\alpha} \}$ where $\sigma (i^\alpha_j )=i^\alpha_{j-1}$ for $1< j\le n_\alpha$. Then
  (see \cite{Lo} 9.2.2)   
  $$\theta_\sigma (\otimes f_i   )=\Pi_\alpha   \tr (f_{i^\alpha_1}\ldots f_{i^\alpha_{n_\alpha}}). \tag 3.3.2$$ 
  
  The map $\theta_I$ is surjective if $r\ge n$: We need to find
 for every $\sigma \in \Sigma_I$ some     $f^\sigma_i \in \End (V)$, $i\in I$, such that $\theta_I (\otimes f^\sigma_i )=\sigma$. Pick a base $\{ v_a  \}$ of $V$ indexed by elements $a$ of $I\sqcup J$,  and define  $f^\sigma_i $ as  $f^\sigma_i (v_i )=v_{\sigma (i)}$,  $f^\sigma_i (v_{a} )=0$ for $a\neq i$.

The  map $\theta_{RI}:=\theta_I \otimes \id_{(R[1])^{\otimes I}}: \End (V)^{\otimes I}\otimes (R[1])^{\otimes I} \to \Bbb Z [\Sigma_I ]\otimes (R[1])^{\otimes I}$  commutes with the diagonal $\Sigma_I$-actions (where $\Sigma_I$ acts on $(R[1])^{\otimes I}$ transposing the factors). One has $\End (V)\otimes R =\End_R (V_R )$, $V_R :=V\otimes R$, so $\End (V)^{\otimes I}\otimes (R[1])^{\otimes I}=(\End (V)\otimes R[1])^{\otimes I}= (\End_R (V_R )[1])^{\otimes I}$. So for $V=\Bbb Z^r$ we get a $\Sigma_I$-equivariant map $$\theta_{RI}: (\fgl_r (R)[1])^{\otimes I} \to \Bbb Z [\Sigma_I ]\otimes (R[1])^{\otimes I}. \tag 3.3.3$$

\proclaim{Lemma} The coinvariants $(\Bbb Z [\Sigma_I ]\otimes (R[1])^{\otimes I})_{\Sigma_I}$ identify naturally with the degree $n$ component of $\text{\rm Sym}\, ( C^\lambda (R)[1])$.
\endproclaim

\demo{Proof of Lemma} Since $\Sigma_I$ acts on $\Bbb Z [\Sigma_I ]$ by conjugation,  $(\Bbb Z [\Sigma_I ]\otimes (R[1])^{\otimes I})_{\Sigma_I}$ equals the direct sum, labeled by  conjugacy classes of  $\sigma\in \Sigma_I$, of  coinvariants $((R[1])^{\otimes I})_{S_\sigma}$; here $S_\sigma$ is the centralizer of $\sigma$ in $\Sigma_I$. For an orbit $I_\alpha $   of $\sigma$ let $\sigma_\alpha \in S_\sigma$ be equal to $\sigma$ on $I_\alpha$ and identity outside; set $n_\alpha :=|I_\alpha |$. The subgroup $N_\sigma = \Pi_\alpha \,\Bbb Z/n_\alpha$  generated by $\sigma_\alpha$'s is normal in $S_\sigma$. 
 The action of $S_\sigma$ on the set of $I_\alpha$'s yields an identification  $S_\sigma /N_\sigma = \Pi_\ell \,\Sigma_{J_\ell}$; here $J_\ell$ is the set of orbits $I_\alpha$ with $n_\alpha =\ell$.
 One has $((R[1])^{\otimes I})_{N_\sigma}= \otimes_\alpha (C^\lambda_{n_\alpha -1} (R)[n_\alpha ])$, hence
$((R[1])^{\otimes I})_{S_\sigma}=\otimes_\ell  ((C^\lambda_{\ell -1} (R)[\ell ])^{\otimes J_\ell})_{\Sigma_{J_\ell}}$, q.e.d. \qed
\enddemo

Since $((\fgl_r (R)[1])^{\otimes I})_{\Sigma_I}= \Lambda^n (\fgl_r (R))[n]$ and by the lemma, the $\Sigma_I$-coinvariants
map  $(\theta_{RI})_{\Sigma_I}$  can be rewritten as $\kappa^r_{Rn}: C (\fgl_r (R))_n \to (\Sym\, ( C^\lambda  (R)[1]))_n$. One checks using (3.3.2) that $\kappa^r_{R} : C  (\fgl_r (R)) \to \Sym\, ( C^\lambda  (R)[1])$ commutes with the differential; it commutes with the coproducts by construction. The surjectivity assertion follows from surjectivity of $\theta_I$. We are done. \qed \enddemo

\subhead{\rm 3.4}\endsubhead Consider the subcomplex $\Ker (\kappa^r_{R} )$ of $C  (\fgl_r (R))$.

\proclaim{Theorem} For every $r\ge 0$ there is a nonzero integer $c_{r}$  that     kills the homology group $H_n \Ker \, (\kappa^r_{R} )$ for every ring $R$ and $n\le r$.     \endproclaim

\demo{Proof}  As in 3.3, set $V=\Bbb Z^r$. Consider the action of the Lie algebra $\fg :=\fgl (V)=\fgl_r (\Bbb Z)$ on $\End (V)^{\otimes n}$. The map $\theta_n : \End (V)^{\otimes n}\twoheadrightarrow \Bbb Z [\Sigma_n ]$ from 3.3 commutes with the $\fg $-actions, so $\Ker (\theta_n )$ is a free $\Bbb Z$-module with $\fg  $-action.  
By the invariant theory \cite{Lo} 9.2.8 and  reductivity of $\fg_{\Bbb Q}$,   $\Ker (\theta_n )_{\Bbb Q}$ is a direct sum of finitely many {\it nontrivial} irreducible $\fg_{\Bbb Q }$-modules. Let $\fc \in U(\fg )$ be the Casimir element. Recall that $\fc$ lies in the center of the enveloping algebra $U (\fg ) $, it kills the trivial representation, i.e., $\fc \in \fg U (\fg )$, and its action on any nontrivial irreducible finite-dimensional representation of $\fg_{\Bbb Q}$ is nontrivial. So we can find $h(t)\in t\Bbb Q[t]$ such that $h(\fc )$ acts as identity on  $\Ker (\theta_n )_{\Bbb Q}$. Let $c$ be an integer such that $cg(t)\in t\Bbb Z[t]$. 

We check that $(n! c)^2$ kills  $H_n \Ker \, (\kappa^r_{R} )$; then $c_r :=(c^{r}\Pi_{1\le n\le r} n!)^{2}$ satisfies the condition of the theorem. Let
 $\pi : (\fgl_r (R)[1])^{\otimes n}\twoheadrightarrow   ((\fgl_r (R)[1])^{\otimes n})_{\Sigma_n}  =C (\fgl_r (R))_n [n]$ be the projection and $s: 
C (\fgl_r (R))_n [n]\to (\fgl_r (R)[1])^{\otimes n}$ be the map such that $s\pi = \Sigma_{\sigma \in \Sigma_n} \sigma$. Then $\pi$ and $s$ commute with  $\fg$-action,   $\pi s$ is multiplication by $n!$.   Since $\theta_n$ is surjective, one has
 $s(\Ker \, (\kappa^r_{Rn} ))\subset \Ker (\theta_{Rn})=\Ker (\theta_n )\otimes (R[1])^{\otimes n}$. Thus $h(\fc )s $ equals $cs$ on $\Ker \, (\kappa^r_{Rn} )$; composing with $\pi$, we see that $n!h(\fc )$ equals $n! c$ on $\Ker \, (\kappa^r_{Rn} )$. The adjoint action of $\fgl_r (R)$ on the Chevalley complex $C  (\fgl_r (R))$ is homotopically trivial, so such is the action of the subalgebra $\fg$. So, since 
$h(\fc )\in\fg U(\fg )$, its action on $C  (\fgl_r (R))$ is homotopic to zero; since $h(\fc ) $ sends 
$C  (\fgl_r (R))$ to $\Ker \, (\kappa^r_{R} )$, the action of $h(\fc )^2 $ on
$\Ker \, (\kappa^r_{R} )$ is homotopic to zero. The multiplication by $(n!c)^2$ on $\Ker \, (\kappa^r_{R} )_n$ equals $  (n!)^2 h(c)^2$, hence it  kills 
the homology. \qed
\enddemo

 \subhead{\rm 3.5}\endsubhead For a $p$-adic Lie algebra $\fg=\limleft \fg /p^i \fg$  we have the Chevalley pro-complex $C(\fg\\hat\, ):=\qlimleft C(\fg /p^i \fg )$, etc. So for a
 $p$-adic $A$ as in 3.2 we have pro-complexes
  $C(\fgl_r (A)\\hat \, )  :=   \qlimleft  C  (\fgl_r (A_i ))$, 
 $ \Sym  (C^\lambda   (A)\\hat \, [1]) :=\qlimleft \Sym  (C^\lambda   (A_i )  [1])$, and a morphism  $\kappa^r_{A\\hat}=\qlimleft  \kappa^r_{A_i}: C (\fgl_r (A)\\hat \, )\to \Sym  (C^\lambda   (A)\\hat \, [1])$.

\proclaim{Corollary} The map $\tau_{\le r}\kappa^r_{A\\hat} :  \tau_{\le r} C (\fgl_r (A)\\hat \, )  \to \tau_{\le r}  \text {\rm Sym}  (C^\lambda   (A)\\hat \, [1])$ is a quasi-isogeny. \qed
\endproclaim

\remark{Remark}  Suppose $A$ has bounded $p$-torsion.  Then the terms of 
$C (\fgl_r (A)\\hat \, )$ and $\Sym  (C^\lambda   (A)\\hat \, [1])$ have bounded $p$-torsion as well, so
 the evident maps $C (\fgl_r (A)\\hat \, )\otimes^L C (\fgl_r (A)\\hat \, )\to C (\fgl_r (A)\\hat \, )\otimes C (\fgl_r (A)\\hat \, )$, $\Sym  (C^\lambda   (A)\\hat \, [1])\otimes^L \Sym  (C^\lambda   (A)\\hat \, [1])\to \Sym  (C^\lambda   (A)\\hat \, [1])\otimes \Sym  (C^\lambda   (A)\\hat \, [1])$ are isogenies. So the coproducts on  $C (\fgl_r (A)\\hat \, )$ and $\Sym  (C^\lambda   (A)\\hat \, [1])$ yield maps 
$\Delta_C : C (\fgl_r (A)\\hat \, )_{\Bbb Q} \to C (\fgl_r (A)\\hat \, )_{\Bbb Q}\otimes^L C (\fgl_r (A)\\hat \, )_{\Bbb Q}$, $\Delta_S : \Sym  (C^\lambda   (A)\\hat \, [1])_{\Bbb Q}\to 
\Sym  (C^\lambda   (A)\\hat \, [1])_{\Bbb Q}\otimes^L \Sym  (C^\lambda   (A)\\hat \, [1])_{\Bbb Q}$ in pro-$\CD (\CA b )^-_{\Bbb Q}$, and one has
$(\kappa^{r}_{A\\hat})^{\otimes 2}\Delta_C = \Delta_S \kappa^r_{A\\hat}\,$.
\endremark

 \subhead{\rm 3.6}\endsubhead We want to pass to the limit $r\to\infty$ in the above corollary to get rid of the truncation $\tau_{\le r}$. Since the torsion exponent  $c_r$ of the theorem in 3.4 depends on $r$, we cannot pass to the limit directly, and   proceed as follows:
 
 In the setting of 3.3  consider the standard embeddings $\fgl_1 (R) \subset \fgl_2 (R) \subset\ldots$; set $\fgl (R):=\cup \fgl_r (R)$. Maps $\kappa^r_R$   are mutually compatible, so they yield a  morphism  $\kappa_R  :   C  (\fgl  (R)) \to  \Sym ( C^\lambda  (R)[1]).$
The complexes $C  (\fgl_r (R))$ form an increasing filtration on $C  (\fgl (R))$. Denote by $C  (\fgl (R))_{(*)}$ the d\'ecalage of that filtration, so  the $n$th component of  $C (\fgl (R))_{(a)} $ equals $\{ c\in C (\fgl_{a+n} (R))_n : \partial (c)\in C(\fgl_{a+n-1} (R))_{n-1} \}$. Let $\kappa_{aR}$ be the restriction of $\kappa_R$ to $C (\fgl (R))_{(a)}$.

\proclaim{Corollary} For $a\ge 0$ the map $\kappa_{aR}: C (\fgl (R))_{(a)}\to \Sym ( C^\lambda  (R)[1])$ is surjective and the groups  $H_n \Ker (\kappa_{aR})$ are killed by   nonzero integers that depend on $a$ and $n$ but not on $R$. \qed
\endproclaim

Notice that $C (\fgl (R))_{(a)}$ are {\it not} subcoalgebras of $C (\fgl (R))$.

 \subhead{\rm 3.7}\endsubhead
For a $p$-adic $A$ set $C(\fgl (A)\\hat \, )_{(a)} :=   \qlimleft  C  (\fgl (A_i ))_{(a)}$. We have the maps 
 $$\kappa_{aA\\hat }:= \qlimleft  \kappa_{a A_i}  : C  (\fgl (A)\\hat \, )_{(a)} \to  \text{\rm Sym} ( C^\lambda (A)\\hat \, [1]) \tag 3.7.1$$ that are compatible with the embeddings $C  (\fgl (A)\\hat \, )_{(a)} \hra C  (\fgl (A)\\hat \, )_{(a+1)} $.

\proclaim{Corollary} Maps (3.7.1) are  quasi-isogenies for $a\ge 0$. \qed
\endproclaim

\subhead{\rm 3.8}\endsubhead Consider the adjoint action of $GL_r (R)$ on $\fgl_r (R)$ and the corresponding $GL_r (R)$-action on $C(\fgl_r (R))$.

\proclaim{Corollary} For every $r\ge 0$ there is a nonzero integer $d_r$ such that for each ring $R$ and $n\le r$ the   $GL_r (R)$-action  on the subgroup $d_r H_n C(\fgl_r (R))$ of $H_n C(\fgl_r (R))$ is trivial.
\endproclaim

\demo{Proof} Consider the standard embedding $\fgl_r (R)\times \fgl_r (R)\hra \fgl_{2r}(R)$ and the corresponding two embeddings $i_1 ,i_2 :\fgl_r (R)\hra \fgl_{2r}(R)$. Consider the adjoint action of $GL_r (R)$ on $\fgl_{2r}(R)$ via the standard embedding $GL_r (R)\hra GL_{2r}(R)$. The restriction of this action on the image of $i_1$ is the adjoint action on $\fgl_r (R)$, and on the image of $i_2$ it is trivial. The theorem in 3.4 implies that the kernel and cokernel of both maps $i_1 ,i_2 : H_n C(\fgl_r (R))\to H_n C(\fgl_{2r}(R))$ are killed by  some nonzero integer  $b_r$, so we can take $b_r^2$ for $d_r $. \qed
\enddemo

\enspace
\head 4. The Lazard isomorphism  \endhead

 In chapter V of \cite{L} (see \cite{HKN} for a digest) Lazard identifies the continuous group and Lie algebra cohomology for saturated groups of finite rank.
We adapt his argument  to the setting of
 pro-complexes where the finite rank condition becomes irrelevant.

\subhead{\rm 4.1}\endsubhead Let $G$ be a topological group such that
 open normal subgroups $\{ G_\alpha \}$ of $G$ form a base of the topology and the topology is complete, i.e.,  $G=\limleft G/G_\alpha$. We denote  by $C(G \\hat \, )$ the group chain pro-complex $ \qlimleft C(G/G_\alpha ,\Bbb Z)$. The diagonal embedding $G\hra G\times G$ yields a cocommutative coproduct on $C(G \\hat \, )$.

The example we need:
Let $A=\limleft A_i$,   $A_i :=A/p^i$, be a $p$-adic unital ring. Then $GL_r (A)$  carries a  topology as above whose base are congruence subgroups 
 $GL_r (A)^{(i)} :=  \Ker (GL_r (A)\twoheadrightarrow GL_r (A_i ))$. Thus for any open subgroup $G$  of $GL_r (A)$ we have the pro-complex $C(G \\hat \, )$.  As in 3.5, we have the Chevalley pro-complex 
$C(\fgl_r (A)\\hat\,):= \qlimleft  C  (\fgl_r (A_i ))$.

 \proclaim{Theorem} Suppose $A$ has bounded $p$-torsion.  Then for any open $G$ that lies in $ GL_r (A)^{(1)}$ if $p$ is odd and in  $ GL_r (A)^{(2)}$ if $p=2$  there is a natural quasi-isogeny  $$\zeta =\zeta_G   :  C(G\\hat\, )_{\Bbb Q}\iso C(\fgl_r (A)\\hat\,)_{\Bbb Q}.\tag 4.1.1$$ These   $\zeta_G $ are compatible with the coproduct, embeddings of $G$, the adjoint action of $GL_r (A)$,   the embeddings $GL_r (A)\hra GL_{r+1}(A)$, and the morphisms of $A$. 
 \endproclaim

The proof  takes 4.2--4.10. We start with general remarks about the  homology of topological groups of relevant kind   (4.2--4.4), then we recast the theorem, following Drinfeld's suggestion, as a general assertion about   $p$-adic Lie algebras  (4.5--4.6).
The map $\zeta$ is constructed in 4.7, and
we check that it is a quasi-isogeny in 4.8--4.10.

\subhead{\rm 4.2}\endsubhead  Let $K$ be a discrete group. Its finite filtration
$K=K_n \supset   \ldots\supset K_1  \supset K_0 =\{ 1\}$ by normal subgroups is called {\it $p$-special} if every successive quotient $K_i /K_{i-1}$ is an abelian group killed by   $p$. 
We say that $K$ is {\it $p$-special} if it admits a $p$-special filtration. 

\proclaim{Lemma} Suppose that $K$ admits an $n$-step $p$-special filtration. Let $F$ be an abelian group viewed as a $K$-module with trivial $K$-action. Then $H_a (K, F )$, $a>0$, is killed by  $f(a,n)=p^{a+n-1 \choose n-1}$.
\endproclaim

\demo{Proof } Since $C(K,F)=C(K,\Bbb Z)\otimes^L_{\Bbb Z}F$, the group
 $H_a (K,F)$ is isomorphic to a direct sum of $H_a (K,\Bbb Z )\otimes F$ and  Tor$^{\Bbb Z}_1 (H_{a-1}(K,\Bbb Z),F)$. So it is enough to consider the case $F=\Bbb Z$.

 Let $K':=K_{n-1}$ be the next to the top term of the filtration, so  the lemma for $K'$ is known by induction. Let us compute $H_a (K ,\Bbb Z )$ using 
 Hochschild-Serre spectral sequence $E^2_{i,j}=H_i (K/K', H_j (K',\Bbb Z))$. Term $E^2_{i,j}$ is killed by $p^{j+n-2 \choose n-2}$ if $i+j>0$: indeed, for $j>0$ this follows by the induction assumption, and $E^2_{i,0}=H_i (K/K' ,\Bbb Z)$ is killed by $p$ for $i>0$ since $K/K'$ is a direct sum of copies of $\Bbb Z/p$ (use K\"unneth formula and the standard computation of $H_\cdot (\Bbb Z/p ,\Bbb Z )$). Thus $H_a (K ,\Bbb Z )$, $a>0$, is killed by $p^{\Sigma_{0\le j\le a}{ j+n-2 \choose n-2}} = p^{a+n-1 \choose n-1}$, q.e.d. \qed
 \enddemo

\remark{Remarks} (i) If $K' $ is a normal subgroup of $K$, then $K$ is $p$-special if and only if   both $K'$ and $K/K'$ are $p$-special.
\newline (ii) The assertion of the lemma need not be true if $F$ is an arbitrary $K$-module, as follows from the next exercise: \newline {\it Exercise.} For any finite  group $K$ and any $a>0$  find a $K$-module $F$ such that $H_a (K,F)=\Bbb Z/|K|\Bbb Z$. \endremark

\subhead{\rm 4.3}\endsubhead We say that a topological group $G=\limleft G/G_\alpha$ as in 4.1 is {\it $p$-special} if every $G/G_\alpha$ is a $p$-special discrete group. 

\remark{Example}  $GL_r (A)^{(1)}$, hence every its open subgroup, is $p$-special.
\endremark

\remark{Remark} If $G'$ is a closed normal subgroup of $G$ then $G$ 
 is $p$-special if and only if both $G'$ and $G/G'$ are $p$-special.
\endremark

\proclaim{Corollary}   Suppose $G$ is $p$-special. If the map $C(G'\\hat\,)\to C(G\\hat\, )$ is a quasi-isogeny for 
 all $G'$ in some base of open normal subgroups of $G$, then the same is true for every open $G' \subset G$.
\endproclaim

\demo{Proof} Pick a normal open $G'' \subset G'$ such that $C(G''\\hat\,)\to C(G\\hat\, )$ is a quasi-isogeny. Our $G'$ is an extension of $K:=G'/G''$ by $G''$, and we   compute $H_\cdot C (G'\\hat\, )$ using the Hochschild-Serre spectral sequence. It is enough to check that $H_a (K, H_b C(G''\\hat\, ))=\qlimleft H_a (K, H_b (G''/G''\cap G_\alpha ,\Bbb Z ))$ is killed by   a power of $p$ if $a>0$ and that the map $H_b  C(G''\\hat\, ) \to H_0 (K, H_b  C(G''\\hat\, ))$ is an isogeny. We can replace $H_b C(G''\\hat\, )$ by $H_b C(G\\hat\, )$ since  the map $H_b  C(G''\\hat\, ) \to H_b C(G\\hat\, )$ is an isogeny.
Now the action of $K$ on $H_b (G/G_\alpha ,\Bbb Z )$ is trivial, and we are done by the lemma in 4.2. \qed
\enddemo

\subhead{\rm 4.4}\endsubhead It will be convenient to replace integral chains by $p$-adic ones: 
For $G$ as in 4.1 set $C(G)\\hat\, := \qlimleft\!\!_{\alpha ,i} C(G/G_\alpha ,\Bbb Z /p^i )$. There is an evident map  $C(G\\hat\, )\to C(G)\\hat\, $.

\proclaim{Corollary} If $G$ is $p$-special then $\tau_{>0}C(G\\hat\, )\iso \tau_{>0}C(G)\\hat\, $.
\endproclaim

\demo{Proof} Replacing $\Bbb Z/p^i$ by $\CC one\, (p^i :\Bbb Z \to \Bbb Z )$, we see that the homotopy fiber of the arrow is a pro-complex  $\qlimleft\!\!{}_{\alpha ,i} \, \tau_{>0}C (G/G_\alpha ,\Bbb Z_i )$ where $\Bbb Z_i$ are copies of $\Bbb Z$ arranged into a projective system $\ldots\buildrel{p}\over\to \Bbb Z \buildrel{p}\over\to \Bbb Z$. We want to show that its homology  $\qlimleft\!\!{}_{\alpha ,i}\, H_a  (G/G_\alpha ,\Bbb Z_i ) = \qlimleft\!\!_{\alpha}\qlimleft\!\!_{i}\, H_a  (G/ G_\alpha ,\Bbb Z_i ) $, $a>0$, vanish. It is enough to check that $\qlimleft\!\!_{i}\, H_a  (K ,\Bbb Z_i )$ vanishes for every $K=G/G_\alpha$ and $a>0$. By the lemma in 4.2,  $H_a  (K,\Bbb Z )$, $a>0$, is killed by some power of $p$, and we are done by Remark (ii) in 1.2.  \qed
\enddemo

\subhead{\rm 4.5}\endsubhead Let us recast the theorem in 4.1 as a general result about $p$-adic Lie algebras:

Let $\fg$ be any $p$-adic Lie algebra, so $\fg = \limleft \fg/p^i \fg$. Suppose $\fg$ has no $p$-torsion and its Lie bracket 
is divisible by $p$ if $p\neq 2$ and by 4 if $p=2$. Then the  Campbell-Hausdorff series  converges and defines a topological group structure on $\fg$. Denote
the corresponding group by $G_{\fg}$, and let $$\log : G_{\fg} \rightleftarrows \fg : \exp \tag 4.5.1$$ be the mutually inverse ``identity" identifications.  Consider the pro-complexes   $C( G\hat{_\fg} )$, $C(\fg\\hat\, )$ of group and Lie algebra chains (see 4.1, 3.5).

\proclaim{Theorem} There is a canonical quasi-isogeny $$\zeta =\zeta_{\fg} : C(G\hat{_\fg} )_{\Bbb Q} \iso C(\fg \\hat\, )_{\Bbb Q} \tag 4.5.2$$ functorial with respect to morphisms of $\fg$'s. 
\endproclaim

\remark{Remarks} (i) The topological group $G_{\fg}$ is $p$-special. More precisely, set  $\fg^{(i)}:= p^{i-m} \fg $ where $m=1$,
 $i\ge 1$ if $p$ is odd, and $m=2$,
  $i\ge 2$ if $p=2$ (the reason for that normalization will be clear later). Then 
$ G_{\fg}^{(i)}:= \exp ( \fg^{(i)} )$ are normal subgroups that form a basis of the topology on $G_{\fg}$, the successive quotients $\gr^{(\cdot )} G_{\fg}$ are abelian, and one has a canonical identification of abelian groups
$ \gr^{(\cdot )}\exp :    \gr^{(\cdot )} \fg  \iso \gr^{(\cdot )} G_{\fg}$. \newline (ii)  Since $\exp (pa)=\exp(a)^p$, the map $g\mapsto g^p$ yields bijections
$G^{(i)}_{\fg} \iso G^{(i+1)}_{\fg}$ and $\Bbb F_p$-vector space isomorphisms
$\gr^{(i )} G_{\fg}\iso \gr^{(i+1)} G_{\fg}$.
\newline (iii) Consider      the $p$-adic Chevalley  pro-complex  $ C(\fg )\\hat \, :=\qlimleft C(\fg )\otimes \Bbb Z/p^i =  \qlimleft C(\fg/p^i \fg ,\Bbb Z /p^i )$. One has $\tau_{>0}C( \fg \\hat\, )=\tau_{>0}C( \fg )\\hat\, $. Since  the chain complexes in (4.5.2) are direct sums of $\Bbb Z$ and their $\tau_{>0}$ truncations, the corollary in 4.4 implies that (4.5.2) amounts to a natural quasi-isogeny $$\zeta  =\zeta _{\fg}  :   C(G_{\fg})\\hat_{\Bbb Q}\iso C( \fg )\\hat_{\Bbb Q} . \tag 4.5.3$$
\endremark

\remark{Question {\rm (Drinfeld)}} Is there an assertion for $p$-torsion Lie algebras behind the theorem? More precisely, let $\fg$ be a $p$-torsion Lie algebra, and suppose we have  another Lie bracket on $\fg$ such that the original bracket equals $p$ times the new one if $p$ is odd or 4 times the new one if $p=2$. Then we have the group $G_{\fg}$ as above, and there is a natural map $C(G_{\fg})\to C(\fg' )$, where  $\fg'$ is $\fg$ with the new Lie bracket, defined as in 4.7 below. What can one say about this map? 
\endremark

\proclaim{{\rm 4.6.} Lemma}  The  theorem in 4.5 implies that in 4.1.
\endproclaim

\demo{Proof} (a) It is enough to consider the case when $A$
 has no torsion: Indeed, for an arbitrary $A$ let $A_{\text{tors}} \subset A$ be the ideal of $p$-torsion elements,  $\bar{A}:=   A/A_{\text{tors}}$, and let $\bar{G}$ be the image of $G$ in $GL_r (\bar{A})$. Then  $G$ is an extension of $\bar{G}$ by $K:=G\cap (1+ \text{Mat}_r (A_{\text{tors}} ))$, and $K$ is discrete since $A$ has bounded $p$-torsion. Since  $K$ is $p$-special (see the remark in 4.3), the Hochschild-Serre spectral sequence and 4.2 show that the map $C(G\\hat\, )\to C(\bar{G}\\hat\, )$ is a quasi-isogeny. The map $C(\fgl_r (A)\\hat\, )\to C(\fgl_r (\bar{A})\\hat \, )$ is a quasi-isogeny as well, and we are done.
 \newline (b) It is enough to construct quasi-isogeny $\zeta_G$ of (4.1.1) when $G$ is a congruence subgroup: Indeed, by the corollary in 4.3, we know then that
embeddings of $G$'s yield quasi-isogenies between the pro-complexes $C(G\\hat\, )$, which gives $\zeta_G$ for every $G$. \newline(c)  Suppose that $m\ge 1$ if $p\neq 2$ and $m\ge 2$ if $p=2$. The logarithm and exponential series converge $p$-adically on $GL_r (A)^{(m)}$ and $p^m \fgl_r (A)$ and define mutually inverse continuous bijections $$\log : GL_r (A)^{(m)} \rightleftarrows p^m \fgl_r (A) : \exp  \tag 4.6.1$$ that satisfy the Campbell-Hausdorff formula. One can interpret this saying that for $\fg =  p^m \fgl_r (A)$ the corresponding $G_{\fg}$ identifies canonically with $GL_r (A)^{(m)}$ preserving the logarithm and exponent maps. 

So (4.5.2) for $\fg = p^m \fgl_r (A)$ can be rewritten as a quasi-isogeny $C( GL_r (A)^{(m)}\\hat\,)_{\Bbb Q} $ $ \iso C (p^m \fgl_r (A)\\hat\, )_{\Bbb Q}$. Now the map $C (p^m \fgl_r (A)\\hat\, )\to C ( \fgl_r (A)\\hat\, )$ that comes from the embedding $p^m \fgl_r (A)\hra \fgl_r (A)$ is a quasi-isogeny as well (which is clear since   $  C (p^m \fgl_r (A)\\hat\, )_n = p^{mn} 
C ( \fgl_r (A)\\hat\, )_n$). We define
 $\zeta_{ GL_r (A)^{(m)}}$ from (4.1.1) as  composition of the above quasi-isogenies, and we are done. The compatibility assertions in the theorem 4.1 follow since $\zeta_{\fg}$ is functorial with respect to morphisms of $\fg$'s  \qed \enddemo

\subhead{\rm 4.7}\endsubhead {\it Proof of the theorem in 4.5.} From now on $\fg$ is a Lie algebra as in 4.5 and $G=G_{\fg}$. Let us
 construct  the map $\zeta_{\fg} : G(G)\\hat_{\Bbb Q} \to C(\fg )\\hat_{\Bbb Q}$ of (4.5.3). 

Let $U(\fg )\\hat\, :=\limleft U(\fg )_i$, where $U(\fg )_i = U(\fg )/p^i U(\fg )=  U_{\Bbb Z/p^i }(\fg /p^i \fg )$, be the $p$-adic completion of the enveloping algebra $U(\fg )$. One has
$$C(\fg )\\hat = \Bbb Z_p \hotimes^L_{U(\fg )\\hat \,}\Bbb Z_p := \qlimleft (\Bbb Z /p^i )\otimes^L_{U(\fg )_i}(\Bbb Z/p^i ).\tag 4.7.1$$

Set $\fg' :=p^{-m}\fg\subset \fg\otimes\Bbb Q$ where $m=1$ if $p\neq 2$ and $m=2$ if $p=2$. The conditions on the Lie bracket on $\fg$ mean that it extends to the Lie bracket on $\fg'  \supset \fg$. We have same objects as above for  $\fg'$, and
$$C(\fg' )\\hat = \Bbb Z_p \hotimes^L_{U(\fg' )\\hat \,}\Bbb Z_p := \qlimleft (\Bbb Z /p^i )\otimes^L_{U(\fg' )_i}(\Bbb Z/p^i ).\tag 4.7.2$$

 Similarly,
consider the Iwasawa algebra $\Bbb Z_p [[G]]:=\limleft (\Bbb Z/p^i )[G/G^{(i)} ]$ where $G^{(i)}$ are as in Remark (i) in 4.5. Then $$C(G)\\hat =  \Bbb Z_p \hotimes^L_{\Bbb Z_p [[G]]}\Bbb Z_p := \qlimleft (\Bbb Z /p^i )\otimes^L_{(\Bbb Z/p^i )[G/G^{(i)} ]}(\Bbb Z/p^i ). \tag 4.7.3$$

One has natural maps of topological algebras $$\Bbb Z_p [[G]] \buildrel{\eta}\over\to U(\fg' )\\hat\,\buildrel{\iota}\over\leftarrow U(\fg )\\hat \tag 4.7.4$$ where $\iota$ comes from the   embedding $\fg\hra\fg'$ and $\eta$ is defined as follows. The exponential series in
$U(\fg')\\hat$ converges on $p^m U(\fg' )\\hat\,$ and yields  a map $$\exp_U : p^m U(\fg' )\\hat\, \to (1+p^m  U(\fg')\\hat\,  )\subset  U(\fg')\\hat\, {}^\times \tag 4.7.5$$ that
satisfies the Campbell-Hausdorff formula. We get a continuous group homomorphism (here  $\log (g) \in \fg \subset \fg'$ is as in  (4.5.1)) $$\eta  : G \to U(\fg')\\hat\, {}^\times, \quad \xi  (g):=\exp_U (\log (g)). \tag 4.7.6$$ 

 We define $\eta$ in (4.7.4) as the  morphism of topological rings that extends (4.7.6). Notice that $\eta  =\limleft \eta_i$ where $\eta_i  : (\Bbb Z/p^i )[G/G^{(i )}]\to U ( \fg' )_i =U_{\Bbb Z/p^i}( \fg'/p^i \fg'  )$. 

Applying (4.7.4)  to  the Tor pro-complexes of (4.7.1)--(4.7.3) we get the maps $$G(G)\\hat \buildrel{\xi}\over\to C(\fg' )\\hat\, \buildrel{i}\over\leftarrow C(\fg )\\hat . \tag 4.7.7$$

The map $i : C(\fg )\\hat \, \to C(\fg' )\\hat\,$ is evidently a quasi-isogeny since it identifies
$C(\fg )\\hat_n$ with $p^{mn} C(\fg' )\\hat_n$. We define $\zeta_{\fg}   : G(G)\\hat_{\Bbb Q} \to C(\fg )\\hat_{\Bbb Q}$ of (4.5.3) as the composition $i^{-1}\xi$ of (4.7.7).

 It remains to prove that the Tor map $\xi : G(G)\\hat \,\to C(\fg' )\\hat \,  $ for $\eta$  is a quasi-isogeny.

\subhead{\rm 4.8}\endsubhead Let $I_p := \Ker (\Bbb Z [G]\to\Bbb Z/p )$ be the $p$-augmentation ideal; set $\Bbb Z_p [[G]]\,\bar{}\, := \limleft \Bbb Z [G]/I_p^i$.
 Remark (ii) in 4.5 implies that  the  $I_p$-adic topology is weaker
than the Iwasawa algebra topology from 4.7,\footnote{By loc.~cit., it is enough to find for any  $n>0$ some $a=a(n)$ such that $g^{p^a}$ equals 1 in $\Bbb Z[G]/ I_p^n$ for every $g\in G$. If $p^b >n$ then $g^{p^b}$ equals 1 in $\Bbb Z[G] / I_p^n + p \Bbb Z[G]$,  so we can take $a= b+n $. } so we have a continuous  map
$\alpha: \Bbb Z_p [[G]]\to \Bbb Z_p [[G]]\,\bar{}\,$, $\alpha (g)=g$. Since $\eta (I_p )\subset p U(\fg' )\\hat\,$, $\eta$ is continuous for the $I_p $-adic topology, i.e., we have a continuous map $\bar{\eta} : \Bbb Z_p [[G]]\,\bar{}\,\to U(\fg')\\hat \,$ such that $\eta= \bar{\eta} \alpha$.
Set $$C(G)\,\bar{} =\Bbb Z_p \hotimes^L_{\Bbb Z_p [[G]]\,\bar{}}\,\, \Bbb Z_p := ``\limleft\! " \, (\Bbb Z /p^i )\otimes^L_{\Bbb Z [G]/I_p^i} (\Bbb Z /p^i ),  \tag 4.8.1$$ and let $\tau$, $\bar{\xi} $ be the Tor maps for $\alpha$ and $  \bar{\eta}$. Since $\xi  = \bar{\xi}  \tau$,  the theorem in 4.5, hence that in  4.1, follows from the next result to be proven in 4.9--4.10 below:

\proclaim{Theorem} (i) The map $\bar{\xi}  : C(G) \,\bar{}\,\to C(\fg' )\\hat$ is a quasi-isogeny. \newline (ii) The map $\tau: C(G)\\hat \to C(G)  \,\bar{}\,$ is a quasi-isomorphism.
\endproclaim

\remark{Remark} The map $\alpha : \Bbb Z_p [[G]]\to \Bbb Z_p [[G]]\,\bar{}\,$ is {\it not} an isomorphism unless $\fg$ is a finitely generated $\Bbb Z_p$-module:  Indeed, $\Bbb Z_p [[G]]$ has discrete quotient $\Bbb F_p [ G/ G^{(2)}]  = \Bbb F_p [\fg/p\fg  ]$ (we assume $p$ is odd), while the corresponding quotient of 
$\Bbb Z_p [[G]]\,\bar{}\,$ is the completion of that group algebra with respect to powers of the augmentation ideal.
\endremark

\subhead{\rm 4.9}\endsubhead {\it Proof of 4.8(i)}. We adapt Lazard's argument from chapter V of \cite{L}:

 Below we identify the graded ring $\gr^\cdot \Bbb Z_p $ for
the $p$-adic filtration  with $\Bbb F_p [t]$, $\deg t=1$, where $t$ is class of $p$ in $p\Bbb Z /p^2 \Bbb Z =\gr^1 \Bbb Z_p$. For every flat $\Bbb Z_p$-algebra $B$ the associated graded ring $\gr^\cdot B$  for the $p$-adic filtration on $B$ identifies in the evident manner with $B_1 \otimes \Bbb F_p [t]=B_1 [t]$ where $B_1 := B/pB$. For example, we have $\gr^\cdot U(\fg' )\\hat \, = U_{\Bbb F_p}(\fh  )[t]=U_{\Bbb F_p [t]} (\fh [t] )$ where $\fh := \fg'/p\fg'$. Recall that $\fg =p^m \fg'$. 
 
Let $R$ be the closure of the image of $\bar{\eta}$ (or $\eta$) in $U(\fg' )\\hat\,$. Let $R^{[n]}:= R\cap p^n U(\fg')\\hat \,$ be the ring filtration on $R$ induced by the $p$-adic filtration on $U(\fg' )\\hat \,$ and $\gr^{[\cdot ]}R\subset \gr^\cdot U(\fg' )\\hat \, $ be the associated graded ring.

 \proclaim{Lemma} (i) The map $\bar{\eta}: \Bbb Z_p [[G]]\,\bar{}\, \to R$ is a homeomorphism of topological rings. \newline (ii) $\gr^{[\cdot ]}R$  is the $\Bbb F_p [t]$-subalgebra of $U_{\Bbb F_p}(\fh )[t]$ generated by $t^m \fh$ which is $U_{\Bbb F_p [t]} (t^m \fh [t] )$ $\subset U_{\Bbb F_p [t]} (\fh [t] )$. \endproclaim

\demo{Proof} Let $I:= \Ker (\Bbb Z [G]\to \Bbb Z )$ be the augmentation ideal, and let $\Bbb Z [G]^{[n]}$ be the sum of the ideals $p^a I^b$ where $a,b\ge 0$ and $a+mb \ge n\ge 0$.
This is a ring filtration on $\Bbb Z [G]$, and the  topology it defines equals the $I_p$-adic topology.\footnote{Indeed, one has $\Bbb Z [G]^{[mn]} \subset I_p^n \subset \Bbb Z [G]^{[n]}$.}  The  ring $\gr^{[\cdot ]}\Bbb Z [G]$ is a graded $\Bbb F_p [t]$-algebra generated by $\CG := I/(I^2 +pI )\subset \gr^{[m]}\Bbb Z [G]$. 

One has $\eta ( \Bbb Z [G]^{[\cdot ]} )\subset R^{[\cdot ]}$, so we have the map of   graded rings  $\gr^{[\cdot ]} \eta :  \gr^{[\cdot ]} \Bbb Z [G] \to   \gr^\cdot U(\fg' )\\hat \,$.  
To prove the lemma it is enough to show that $\gr^{[\cdot ]} \eta $ is injective.
 
Recall that $I/I^2 =G/[G,G]$, so $\CG = G/G^p [G,G]$ that identifies canonically with $\fg /p\fg$ by Remarks (i), (ii) in 4.5. The map $\gr^{[m]}\eta |_\CG$ is  the embedding $\fg/p\fg =t^m \fh \hra U_{\Bbb F_p}(\fh )[t]$. Therefore the $\Bbb F_p [t]$-submodule of $\gr^{[\cdot ]} \Bbb Z [G]$   generated by $\CG$ equals $\CG [t]=\Bbb F_p [t]\otimes\CG$ and  $\gr^{[\cdot ]} \eta $ identifies it with $t^m\fh [t]$.

Notice that $[\CG ,\CG ]$ lies in $t^m \CG \subset  \gr^{[2m]}\Bbb Z [G]$. Thus $\CG [t]$ is a graded Lie $\Bbb F_p [t]$-subalgebra of $\gr^{[\cdot ]} \Bbb Z [G]$,  so we have a surjective map of graded associative algebras $\nu : U_{\Bbb F_p [t]}(\CG [t] )\twoheadrightarrow \gr^{[\cdot ]} \Bbb Z [G]$. The composition $(\gr^{[\cdot ]} \eta )\nu :  U_{\Bbb F_p [t]}(\CG [t] )\to U_{\Bbb F_p}(\fh )[t]$ is the morphism of enveloping algebras that corresponds to the embedding of Lie algebras $\CG [t] \iso t^m \fh [t] \hra \fh [t]$. It is injective, hence $\gr^{[\cdot ]} \eta$ is injective, q.e.d.
    \qed
\enddemo

Set $C(R) \,\bar{}\, := \qlimleft (\Bbb Z/p^i )\otimes^L_{R/R^{[i]}}\! (\Bbb Z /p^i )$, so the embedding $R\hra U(\fg' ) \\hat \,$ yields a
 map of pro-complexes $ C(R)\,\bar{}\,  \to C(\fg' )\\hat \,   $ (see (4.7.2)). 
By (i) of the lemma,   4.8(i) amounts to the next assertion:

\proclaim{Proposition} The morphism  $ C(R)\,\bar{}\,  \to C(\fg' )\\hat \,   $ is a quasi-isogeny.
\endproclaim

\demo{Proof} 
Let $K$ be the Chevalley chain complex of the Lie $\Bbb F_p [t]$-algebra $t^m \fh [t]$ with coefficients in the free module  $U_{\Bbb F_p [t]} (t^m \fh [t])$. This is a graded free  $U_{\Bbb F_p [t]} (t^m \fh  [t] )$-resolution of $\Bbb F_p [t]$ (with trivial Lie algebra action) with terms $K_n = t^{mn} \Lambda^n_{\Bbb F_p} (\fh )\otimes_{\Bbb F_p} U_{\Bbb F_p [t]} (t^m \fh [t] )$.
Our $K$ can be lifted to a {\it complete filtered} $R$-resolution $P$ of $\Bbb Z_p$,\footnote{One constructs  silly truncations of $P$ by induction using $\gr^{[\cdot ]} \Ker (\partial )= \Ker (\gr^{[\cdot ]} \partial )$ where $\partial$ is the differential of $P$.  
} so $$  K \iso \gr^{[\cdot ]} P. \tag 4.9.1$$ Thus $P_n$ is a complete filtered free $R$-module with generators in filtered degree $mn$.

Let $S$ be the induced complete filtered $ U (\fg' ) \\hat \,  $-complex, so $S=\limleft S/S^{[i]}$ where $S/S^{[i]} = U_{\Bbb Z /p^i} (\fg'/p^i\fg' )\otimes_{R} (P/P^{[i]})$. Notice that $S^{[mn+i]}_n = p^i S_n$ for $i\ge 0$. Set $T_n := p^{-mn} S_n \subset \Bbb Q_p \otimes_{\Bbb Z_p}S_n$; then $T$
is a subcomplex of 
$ \Bbb Q_p \otimes_{\Bbb Z_p}S$. We equip $T$ with the $p$-adic filtration, $T^{ i}= p^i T$. Then (4.9.1) yields isomorphism $$ L[t]
 \iso \gr^\cdot T \tag 4.9.2$$ where $L$ is the Chevalley complex of  $\fh$ with coefficients in $U_{\Bbb F_p} (\fh )$.

Using resolutions $P$ and $T$, one
 can realize the homotopy map of the proposition as  composition 
 $C(R)\,\bar{}\, \isol \qlimleft  (\Bbb Z/p^i )\otimes^L_{R/R^{[i]}} (P/P^{[i]} )
 \buildrel\alpha\over\to \qlimleft  (\Bbb Z/p^i )\otimes_{R} (P/P^{[i]} ) =\qlimleft  (\Bbb Z/p^i )\otimes_{U(\fg' )} (S/S^{[i]} )\buildrel\beta\over\to \qlimleft  (\Bbb Z/p^i )\otimes_{U( \fg'  )}( T/T^i )\iso C( \fg' )\\hat\,$.
Here the first quasi-isomorphism $\isol$ comes since $P/P^{[i]}$ is a resolution of $\Bbb Z/p^i$ by (4.9.1), and the last quasi-isomorphism $\iso$ comes since  $\Bbb Z/p^i  \otimes_{U( \fg' )} T/T^i = \Bbb Z/p^i \otimes^L_{U( \fg'/p^i\fg' )}\Bbb Z/p^i = C(\fg'/p^i \fg' ,\Bbb  Z/p^i  )$ by (4.9.2). The map $\beta$ is a quasi-isogeny  by the construction of $T$. 

  To finish the proof, we  show that $\alpha$ is a quasi-isomorphism.
It is enough to check that $\qlimleft \!\!{}_i \, \text{Tor}_a^{R/R^{[i]}}\! (\Bbb Z/p^j , P_n /P_n^{[i]} )=0$   for every $n$, $j$ and any $a>0$. This follows since 
the map $\text{Tor}_a^{R/R^{[mn+i]}}\! (\Bbb Z/p^{j }, P_n /P_n^{[mn+i ]}) \to \text{Tor}_a^{R/R^{[i]}}\! (\Bbb Z/p^j , P_n /P_n^{[i]} )$ factors through $\text{Tor}_a^{R/R^{[i]}}\! (\Bbb Z/p^j , P_n /P_n^{[mn+i]} )$ and  $P_n /P_n^{[mn+i]}$ is a free $R/R^{[i]}$-module.  \qed
\enddemo

\subhead{\rm 4.10}\endsubhead {\it Proof of 4.8(ii).} (a)  Set $C(G,\Bbb Z/p^j )\\hat   := \qlimleft\!\!{}_i \, (\Bbb Z /p^j )\otimes^L_{(\Bbb Z/p^i )[G/G^{(i)}]} (\Bbb Z /p^i )=   \qlimleft\!\!{}_i \, C(G/G^{(i)},\Bbb Z/p^j )$ and $C(G,\Bbb Z/p^j )\,\bar{}\, := \qlimleft\!\!{}_i \, (\Bbb Z /p^j )\otimes^L_{\Bbb Z [G]/I_p^i } (\Bbb Z /p^i )$. Since $C(G)\\hat\, $ and $C(G )\,\bar{}\,$ are homotopy limits of pro-complexes $
C(G,\Bbb Z/p^j )\\hat \,$ and $C(G,\Bbb Z/p^j )\,\bar{}\,$, see  (4.7.3) and (4.8.1),   assertion 4.8(ii) would follow if we show that the natural maps $\tau : C(G,\Bbb Z/p^j )\\hat  \to C(G,\Bbb Z/p^j )\,\bar{}\,$ are quasi-isomorphisms. By devissage, it is enough to check that $\tau : C(G,\Bbb Z/p )\\hat  \to C(G,\Bbb Z/p )\,\bar{}\,$ is a quasi-isomorphism.

(b) We compute $C(G,\Bbb Z/p )\,\bar{}\,$ using  resolution $P$ from the proof of the proposition in 4.9: As in loc.~cit., we have quasi-isomorphisms $C(G,\Bbb Z/p )\,\bar{}\, \isol \qlimleft  (\Bbb Z/p )\otimes^L_{R/R^{[i]}} (P/P^{[i]} )  \buildrel\alpha\over\to \qlimleft  (\Bbb Z/p )\otimes_{R} (P/P^{[i]} )$ where $\alpha$ is a quasi-isomorphism by the argument at the end of 4.9. By the construction,  $\qlimleft  (\Bbb Z/p )\otimes_{R} (P/P^{[i]} )$ equals $(\Bbb Z /p )\otimes_{U_{\Bbb F_p [t]} (t^m \fh [t])}  K$ which is the same as the Chevalley complex of the {\it abelian} Lie $\Bbb F_p$-algebra $t^m \fh [t]/t^{m+1} \fh [t]$  with coefficients in $\Bbb F_p$. Thus, using the multiplication by $t^m$ isomorphism $\fh \iso t^m \fh [t]/t^{m+1} \fh [t]$, we get identification
 $H_n C(G,\Bbb Z/p )\,\bar{}\, =\Lambda^n_{\Bbb F_p} \fh$.

(c) Let us compute $C(G,\Bbb Z/p )\\hat = \qlimleft   C(G/G^{(i)},\Bbb Z/p )$. This is a cocommutative counital $\Bbb F_p$-coalgebra in the usual way.

\proclaim{Proposition} There is a canonical isomorphism of coalgebras
$$H_\cdot C(G,\Bbb Z/p )\\hat \iso
\Lambda^\cdot_{\Bbb F_p}( \fh ).\tag 4.10.1$$
\endproclaim

\demo{Proof} Set $G_{i,n}:= G^{  (n)}/G^{  (n+i)}$. Below we write $H_\cdot (?)$ for $H_\cdot (?,\Bbb Z/p)$.
 
(i) By Remark (ii) in 4.5 the projection $G_{i,n} \twoheadrightarrow G_{1,n} =\fh$, $\exp(p^n a )\mapsto  a\,\text{mod}\,p\fg'$,  yields an isomorphism $\epsilon_{i,n} : H_1 (G_{i,n}  )\iso H_1 (\fh )=\fh$; here we view $\fh$ as mere abelian group. 
 
 (ii) The kernel of the projection $G_{i+1,n} \twoheadrightarrow G_{i,n}$ is $G_{1, n+i}=\fh$, and the adjoint action of $G_{i+1,n}$ on it is trivial. Thus the corresponding Hochschild-Serre spectral sequence converging to
$H_\cdot ( G_{i+1,n}  )$ has $E^2_{p,q}= H_p (G_{i,n}  )\otimes H_q (G_{1,n+i}  )$. By (i) for $G_{i+1,n}$, the component   $\nu_{i,n} : H_2 (G_{i,n}    )\to H_1 (G_{1,n+i}   )=\fh$ of the differential  in $E^2$ is surjective.

(iii)\footnote{The reader may find it more convenient to follow the computation in (iii) and (iv) using the dual language of cohomology (the coalgebra structure turns then into the algebra one).}  Case $p\neq 2$: Then $\Lambda^\cdot_{\Bbb F_p} (\fh )$ is the free 
cocommutative counital graded $\Bbb F_p$-coalgebra cogenerated by a copy of $\fh$ in degree 1. Since $H_\cdot (G_{i,n}   )$ are cocommutative counital graded $\Bbb F_p$-coalgebras, we have compatible morphisms of graded coalgebras $\tilde{\epsilon}_{i,n}: H_\cdot (G_{i,n}   )\to \Lambda^\cdot_{\Bbb F_p} (\fh )$  that equal  $\epsilon_{i,n}$ in degree 1. We define (4.10.1) as   $\limleft \tilde{\epsilon}_{i,n}$. To see that it is an isomorphism, let us compute the homology of $G_{i,n}$.

The free cocommutative counital graded $\Bbb F_p$-coalgebra cogenerated by a copy of $\fh$ in degree 2 equals $\Gamma_{\Bbb F_p}^{\cdot /2} (\fh )$.\footnote{Here $\Gamma^\cdot_{\Bbb F_p} (\fh )$ is the divided powers Hopf $\Bbb F_p$-algebra generated by a copy of $\fh$.} Thus $\nu_{i.n}$ extends to a morphism of graded coalgebras $\tilde{\nu}_{i,n}: H_\cdot (G_{i,n} )\to  \Gamma_{\Bbb F_p}^{\cdot /2} (\fh )$. Let us show that the morphism $$\tilde{\epsilon}_{i,n}\otimes\tilde{\nu}_{i,n}: H_\cdot (G_{i,n} )\to\Lambda^\cdot_{\Bbb F_p} (\fh )\otimes \Gamma_{\Bbb F_p}^{\cdot /2} (\fh) \tag 4.10.2$$ of graded coalgebras is an isomorphism.

For $i=1$ one has $G_{1,n} =  \fh$, and the assertion follows from the standard calculation of $H_\cdot (\Bbb Z/p  )$ combined with the K\"unneth formula.\footnote{Use the fact that for every homomorphism $\iota_a : \Bbb Z/p \to G_{1,n}$, $\iota_a (x)= (1+ap^n )^x$, $a\in \fh$, one has $\nu_{1,n}(\iota_a (\xi))=a$, where $\xi$ is the standard generator of $H_2 (\Bbb Z/p )$ (which follows since the map $x\mapsto x^p$ on $G_{2,n}$ yields an isomorphism $G_{1,n}\iso G_{1,n+1}$).}
 We proceed   by  induction by $i$ using the   spectral sequence from (ii) to pass from $i$ to $i+1$: By the induction assumption and since we  know that (4.10.2) is an isomorphism for $G_{1,n+i}$, one has $E^2_{\cdot,\cdot} = H_\cdot (G_{i,n}) \otimes H_\cdot (G_{1,n+i})= (  \Lambda^\cdot_{\Bbb F_p} (\fh )\otimes \Gamma_{\Bbb F_p}^{\cdot /2} (\fh ))\otimes ( \Lambda^\cdot_{\Bbb F_p} (\fh )\otimes \Gamma_{\Bbb F_p}^{\cdot /2} (\fh ))$. By the construction of $\nu_{i,n}$, the differential in $E^2$ identifies the second copy of $\fh$ (which lies in $H_2 (G_{i,n})$) with the third copy of $\fh$ (which is in $H_1 (G_{1,n+i})$).  Since the composition of $\fh \hra H_2 (G_{1,n+i})\to H_2 (G_{i+1,n})$ is injective (indeed, its composition with $\nu_{i+1,n}$ equals $\id_{\fh}$),  the map from the fourth copy of $\fh$ (which lies in $H_2 (G_{1,n+i})$) to $E^\infty_{0,2}$ is injective. The compatibility of the differentials with the coproduct implies then that $E^3_{p,q}= \Lambda^p (\fh )\otimes \Gamma^{q/2}_{\Bbb F_p} (\fh )=E^\infty_{p,q}$, therefore $ \tilde{\epsilon}_{i+1,n}\otimes\tilde{\nu}_{i+1,n}  $ is an isomorphism.

Now, since $\nu_{i,n}$ vanishes on the image of $H_2 (G_{i+1,n}  )\to H_2 (G_{i,n}  )$, the transition map  $H_\cdot (G_{i+1,n}  )\to H_\cdot (G_{i,n}  )$ rewritten in terms of (4.10.2) kills the copy of $\fh$ in degree 2. Thus, since the transition map is a map of coalgebras, it equals the composition $\Lambda^\cdot_{\Bbb F_p} (\fh )\otimes \Gamma_{\Bbb F_p}^{\cdot /2} (\fh )\twoheadrightarrow \Lambda^\cdot_{\Bbb F_p} (\fh )\hra \Lambda^\cdot_{\Bbb F_p} (\fh )\otimes \Gamma_{\Bbb F_p}^{\cdot /2} (\fh )$. This implies that (4.10.1) is an isomorphism, q.e.d.

(iv) Case $p=2$: The  free 
cocommutative counital graded $\Bbb F_2$-coalgebra cogenerated by a copy of $\fh$ in degree 1 equals $\Gamma^\cdot_{\Bbb F_2}(\fh )$. The morphism of graded coalgebras $\tilde{\epsilon}_{1,m}: H_\cdot (G_{1,n})\to \Gamma^\cdot_{\Bbb F_2}(\fh )$ that extends $\epsilon_{1,n}$  is an isomorphism (by the standard computation of $H_\cdot  (\Bbb Z/2 )$ and K\"unneth formula).

As in (iii), the free 
cocommutative counital graded $\Bbb F_2$-coalgebra cogenerated by a copy of $\fh$ in degree 2 equals $\Gamma^{\cdot/2}_{\Bbb F_2}(\fh )$, and  $\nu_{i,n}$ extends to a morphism of graded coalgebras $\tilde{\nu}_{i,n}: H_\cdot (G_{i,n} )\to  \Gamma_{\Bbb F_2}^{\cdot /2} (\fh )$.
Since $\nu_{1,m}$ vanishes on the image of  $H_2 (G_{2,n})\to H_2 (G_{1,n})$ and $\Lambda^\cdot_{\Bbb F_2}(\fh )\subset \Gamma^\cdot_{\Bbb F_2}(\fh ) $ is the largest graded subcoalgebra of $H_{\cdot} (G_{1,n})$ on which $\nu_{1,n}$ vanishes, we see that   $\tilde{\epsilon}_{i,n}$ takes values in $\Lambda_{\Bbb F_2}(\fh )$ if $i\ge 2$, i.e., $\tilde{\epsilon}_{i,n}$ yields $\bar{\epsilon}_{i,n}: 
H_\cdot (G_{i,n} )\to\Lambda^\cdot_{\Bbb F_2} (\fh )$. We define (4.10.1) as 
$\limleft \bar{\epsilon}_{i,n}$. 

Now for $i\ge 1$ we have the map $\bar{\epsilon}_{i,n}\otimes\tilde{\nu}_{i,n}: H_\cdot (G_{i,n} )\to\Lambda^\cdot_{\Bbb F_2} (\fh )\otimes \Gamma_{\Bbb F_2}^{\cdot /2} (\fh )$ of graded coalgebras.  One checks that it is an isomorphism in the same way as we did it in (iii), and this implies that (4.10.1) is an isomorphism, q.e.d.
\qed\enddemo

(d) To finish the proof of 4.8(ii), hence of the theorems in 4.5 and in  4.1, it remains to notice that $H_\cdot C(G,\Bbb Z/p )\,\bar{}\,$ is  a cocommutative coalgebra in the same way as $H_\cdot C(G,\Bbb Z/p )\\hat\,$ is, and our map $\tau : H_\cdot C(G,\Bbb Z/p )\\hat\, \to H_\cdot C(G,\Bbb Z/p )\,\bar{}\,$ is a morphism of graded coalgebras. We have identified both coalgebras with $\Lambda^\cdot_{\Bbb F_p} (\fh )$ in (b) and (c). Since $\tau$  equals $\id_{\fh}$ in degree 1, it equals $\id_{ \Lambda^\cdot_{\Bbb F_p} (\fh )}$ in all degrees, q.e.d. \qed

\subhead{\rm 4.11}\endsubhead We are in the setting of 4.1. The theorem in 4.1 can be partially extended  to some other subgroups of $GL_r (A)$ due to  the next proposition:  

Let $G'\subset GL_r (A)$ be any open $p$-special subgroup  (see 4.3). E.g.~the latter condition holds if
the image of $G'$ in $GL_r (A_1 )$ consists of upper triangular matrices with 1 on the diagonal. Let $G \subset G' $ be a  small enough congruence subgroup.

\proclaim{Proposition} The map  $H_n C(G\\hat\,)\to H_n C(G'\\hat\,)$ is an isogeny if $n\le r$.
\endproclaim

\demo{Proof} Consider the adjoint action of $GL_r (A)$ on  $H_n C(G\\hat\, )$.
Combining the theorem in 4.1 with the corollary in 3.8, we see that for some nonzero integer $c_r$ this action is trivial on   $c_r H_n C(G\\hat\, )$. Now compute $H_n C(G')$ using the  Hochschild-Serre spectral sequence for $G\subset G'$, and use the lemma in 4.2. \qed
\enddemo

\remark{Question} Is it true that the adjoint action of $GL_r (A )$ on $C(\fgl_r (A_i ))$ is trivial? If yes, this can replace the reference to 3.8 in the above proof, and the assumption
$n\le r$ can be dropped.
\endremark

\head 5. Relative $K$-groups and   Suslin's stabilization \endhead

\subhead{\rm 5.1}\endsubhead  We identify the symmetric group $\Sigma_r$ with a subgroup of $GL_r$ in the usual way (the transpositions of standard base vectors). For every $\sigma \in \Sigma_r$ let $U^\sigma \subset GL_r$ be the $\sigma$-conjugate of the subgroup of unipotent upper triangular matrices. For a ring $R$ we have Volodin's simplicial set $X_{r}(R):= \cup_{\sigma \in \Sigma_r }B_{U^\sigma (R) }\subset B_{GL_r (R)}$. 
It is clear that $X_r (R)$ is connected, $\pi_1 (X_r (R))=St_r (R)$. As usually, $GL (R)$ is the union of $GL_1 (R)\subset GL_2 (R) \subset \ldots$, $X(R):=\cup X_r (R)$, etc. 

For a cell complex $P$ let $P\to \tau_{\le n}P$ be the $n$th Postnikov truncation of $P$, so $\pi_i (P)\iso \pi_i (\tau_{\le n}P )$ for $i\le n$ and $\pi_{>n}(\tau_{\le n}P ) =0$. Below ${}^+$ is Quillen's $+$ construction. 
Let $\fr$ be the stable rank of $R$.
The next result is due to Suslin \cite{Su}:

\proclaim{ Theorem} (i) The reduced homology $\tilde{H}_n (X_r (R),\Bbb Z )$ vanish for $r\ge 2n+1$. Therefore $X(R)$ is acyclic.
 \newline
(ii) $X(R)$ identifies naturally with the homotopy fiber of  $B_{GL(R)}\to B_{GL(R)}^+$. \newline
(iii) If  $r\ge \max (2n+1, \fr +n)$, then  $H_n (GL_r (R),\Bbb Z )\iso H_n (GL (R),\Bbb Z)$ and  $\tau_{<n}X_r (R) $ identifies naturally with the homotopy fiber of $B_{GL_r (R)}\to \tau_{\le n}B^+_{GL_r (R)}$. 
\endproclaim

\demo{Proof} (i) is \cite{Su} 7.1. (ii) The composition $X(R)\to B_{GL(R)}\to B_{GL(R)}^+$ factors through $X(R)^+$ which is contractible by (i). We get the map from $X(R)$ to the homotopy fiber of $B_{GL(R)}\to B_{GL(R)}^+$. It is a homotopy equivalence by \cite{Lo} 11.3.6. \newline (iii) The composition of maps  $\tau_{\le n}X_r (R) \to B_{GL_r (R)}\to \tau_{\le n}B_{GL_r (R)}^+  $ factors through $\tau_{\le n}X_r (R)^+ $ which is contractible if $r\ge 2n+1$ by (i), hence comes the map from $\tau_{< n}X_r (R) $ to the homotopy fiber of $B_{GL_r (R)}\to \tau_{\le n}B_{GL_r (R)}^+ $; it is a homotopy equivalence by \cite{Su} 8.1. The rest is \cite{Su} 8.2.
\qed
\enddemo

\subhead{\rm 5.2}\endsubhead
If $J\subset R$ is a two-sided ideal, then the relative Volodin's simplicial set $X_r (R,J) $ is defined as the preimage of $X_{r}(R/J)$ by the projection $B_{GL_r (R)} \to B_{GL_r (R/I)}$, i.e., $X_r (R,J) = \cup_{\sigma \in \Sigma_r }B_{U^\sigma (R,J) }\subset B_{GL_r (R)}$ where $U^\sigma (R,J)$ is the preimage of $U^\sigma (R/J)$ in $GL_r (R)$.  So $X_r (R,J)$ is connected,  $\pi_1 (X_r (R,J))=St_r (R/J)\times_{GL_r (R/J)}GL_r (R)$.

 From now on we assume that $J$ is nilpotent. Then $K_0 (R)\iso K_0 (R/J)$ and $GL_r (R) \twoheadrightarrow GL_r (R/J)$, hence $K_1 (R)\twoheadrightarrow K_1 (R/J)$. Thus the relative $K$-spectrum $K(R,J)$ is connected and $\Omega^\infty K(R,J)$ is the homotopy fiber of $B^+_{GL (R)} \to B^+_{GL(R/J)}$. Let $\fr$ be the stable rank of $R$ or of $R/J$ (they coincide since $J$ is nilpotent).

\proclaim{ Proposition} (i)  $X(R,J)^+ \iso \Omega^\infty K(R,J)$, so $C( X(R,J), \Bbb Z )\iso   C(\Omega^\infty K(R,J),\Bbb Z )$.  
(ii)  $H_n C(X_r (R,J) ,\Bbb Z )\iso H_n C(  X(R,J) ,\Bbb Z)$  for $r\ge \max (2n+1, \fr+n )$.
\endproclaim

\demo{Proof} (i) (see \cite{Lo} 11.3.6) 
The projection $B_{GL(R)}\to B_{GL(R/J)}$ is Kan's fibration. So, by 5.1(ii) applied to $R/J$,  $X(R,J)$ is the homotopy fiber of $B_{GL(R)}\to B_{GL(R/J)}^+$. The map $B_{GL(R)}\to B_{GL(R)}^+$ over $B_{GL(R/J)}^+$ yields the map of the homotopy fibers $X(R,J)\to \Omega^\infty K(R,J)$. By 5.1(ii), the homotopy fiber of the latter map equals $X(R)$. Thus, by 5.1(i), $ C(X(R,J), \Bbb Z)\iso   C(\Omega^\infty K(R,J),\Bbb Z)$, hence the assertion. \newline (ii) Repeat the argument of (i) with $B_{GL(?)}^+$ replaced by $\tau_{\le n}B_{GL_r (?)}^+ $ and $X(R,J)$ replaced by $\tau_{\le n}X_r (R,J) $. \qed\enddemo

\remark{Question} Can it be that (ii)  actually holds for  $r>n$?\footnote{If $R$ is a $\Bbb Q$-algebra, the positive answer follows from Goodwillie's theorem \cite{G}, \cite{Lo} 11.3 combined with the Loday-Quillen stabilization \cite{Lo} 10.3.2.} If yes, this would eliminate the assumption of finiteness of stable rank from the corollary in 5.4, hence from the theorem in 2.2.
\endremark

\subhead{\rm 5.3}\endsubhead Suppose $A$ and $I$ are as  in 2.2 and $A$ has bounded $p$-torsion. Choose $m$ such that $p^m \in I$. 
Then $GL_r (A_i )^{(m)} := \Ker (GL_r (A_i )\to GL_r (A_{\min (i,m)}))  \subset U^\sigma (A_i , I_i )$ for any $\sigma$, hence $B_{GL_r (A_i )^{(m)}}\subset X_r (A_i ,I_i )\subset B_{GL_r (A_i )}$. Consider the maps of chain complexes $C(GL_r (A_i )^{(m)} ,\Bbb Z )\to  C(X_r (A_i ,I_i ),\Bbb Z )$. Applying $\qlimleft$, we get a map of pro-complexes  $$C(GL_r (A)^{(m)}\\hat \,  )\to  C(X_r (A,I)\\hat   \, ,\Bbb Z ). \tag 5.3.1$$ 
\proclaim{Proposition} For $n\le r$ the map $H_n C(GL_r (A)^{(m)}\\hat \, ) \to H_n C(X_r (A,I)\\hat \, ,\Bbb Z )$ is an
isogeny.
\endproclaim

\demo{Proof} $X_r (A_i ,I_i )$ is covered by simplicial subsets $B_{U^\sigma (A_i ,I_i )}$. Let $B_{U (A_i ,I_i )}$ be the intersection of several of these subsets. One has $GL_r (A_i )^{(m)} \subset U (A_i ,I_i )$, so $G:=GL_r (A)^{(m)} \subset G':=\limleft 
U (A_i ,I_i )$. The proposition  follows since the maps $H_n C(G\\hat \, )\to H_n C(G'\\hat \, )$ for $n\le r$ are isogenies by 4.11. \qed
\enddemo

\subhead{\rm 5.4}\endsubhead For a ring $R$ consider the standard embeddings $GL_1 (R )\subset GL_2 (R)\subset \ldots$. The complexes $C( GL_r (R ),\Bbb Z )  $ form an increasing filtration on $C( GL (R ),\Bbb Z )$. As in 3.6, denote by $C(GL (R ),\Bbb Z )_{(*)}$ the d\'ecalage of that filtration, so the $n$th component of $C(GL (R ),\Bbb Z )_{(a)}$ equals $\{ c\in C(GL_{a+n} (R),\Bbb Z )_n : \partial (c)\in 
C(GL_{a+n-1} (R),\Bbb Z )_{n-1}\}$. Replacing $GL_r $ by its congruence subgroup we get projective systems of complexes $C(GL (A_i )^{(m)},\Bbb Z)_{(a )}$. Set $C(GL (A)^{(m)}\\hat\,  )_{(a)}:=\qlimleft C(GL (A_i )^{(m)},\Bbb Z)_{(a)}$.

  Maps (5.3.1) for different $r$'s are compatible, so, by 5.2(i), they produce  maps  $$C(GL (A )^{(m)} \\hat \,)_{(a)}\to C(\Omega^\infty K(A,I)\\hat \, ,\Bbb Z ) \tag 5.4.1$$  that are compatible with the  embeddings  $C(GL (A )^{(m)} \\hat \,)_{(a)} \hra C(GL (A )^{(m)} \\hat \,)_{(a+1)}$ and $C(GL (A )^{(m+1)} \\hat \,)_{(a)}\hra C(GL (A )^{(m)} \\hat \,)_{(a)}$. Notice that  the theorem in 4.1  yields quasi-isogenies (see 3.7) $$C(GL (A )^{(m)} \\hat \,)_{(a)\,\Bbb Q} \iso C(\fgl (A)\\hat \, )_{(a)\,\Bbb Q}. \tag 5.4.2$$

\proclaim{Corollary} Suppose the stable rank of $A_1$ is finite. Then  maps  (5.4.1) are  quasi-isogenies for $a\ge 0$.
\endproclaim

\demo{Proof} The stable rank of $A_1$ equals that of $A/I$, so, by 5.3 and 5.2(ii), we know that for given $n$ the map 
$H_n C(GL_r (A)^{(m)} \\hat \, )\to H_n C(\Omega^\infty K(A,I)\\hat  \, ,\Bbb Z)$   is an isogeny if $r$ is large enough.  Now (5.4.2) together with the corollary in 3.7 implies that $C (GL (A)^{(m)} \\hat \, )_{(0)} \to C (GL (A)^{(m)}\\hat \,   )_{(1)} \to \ldots$ are all quasi-isogenies, q.e.d.   \qed\enddemo 

\head 6. Coda
\endhead

It remains to tie up the segments of the tail:

\subhead{\rm 6.1}\endsubhead 
Let $A$ and $I$ be  as in the theorem in 2.2. Denote by $\chi$  the composition of the chain of quasi-isogenies $
 \Sym (CC(A\\hat \, )[1])_{\Bbb Q} \buildrel{3.2}\over\lra \Sym (C^\lambda (A\\hat \, )[1])_{\Bbb Q} \buildrel{(3.7.1)}\over\longleftarrow C(\fgl (A)\\hat \, )_{(a )\,\Bbb Q}$ $ \buildrel{(5.4.2)}\over\longleftarrow C(GL(A)^{(m)}\\hat \,   )_{(a )\,\Bbb Q } \buildrel{(5.4.1)}\over\lra
 C(\Omega^\infty K(A,I)\\hat ,\Bbb Z )_{\Bbb Q} $; it 
 is evidently independent of the choice of $a\ge 0$, $m\ge 1$ involved.
  Consider the composition 
 $$     CC(A\\hat \, )[1]_{\Bbb Q} \hra  \Sym (CC(A\\hat \, )[1])_{\Bbb Q} \iso   C(\Omega^\infty  K(
 A,I)\\hat  ,\Bbb Z )_{\Bbb Q} \to   K( A,I)\\hat_{\Bbb Q} \tag 6.1.1$$ where the first arrow is the evident embedding, the second arrow is $\chi$, and the third one is
 $\nu_{K(A,I)\\hat}$ from  1.2(c).
 
 \proclaim{Proposition} The  composition  $    CC(A\\hat \, )[1]_{\Bbb Q}\to   K( A,I)\\hat_{\Bbb Q} $      is a quasi-isogeny. \endproclaim 

 The promised quasi-isogeny (2.2.1) is its inverse.

 \demo{Proof} Let $\psi :  CC(A\\hat \, )[1]_{\Bbb Q} \to   C(\Omega^\infty  K(
 A,I)\\hat  ,\Bbb Z )_{\Bbb Q}$ be the composition of the first two arrows in (6.1.1).
 According to the proposition in 1.2(c), it is enough to check that $H_n (\psi )$ identifies $H_n (CC  (A\\hat \, )[1]_{\Bbb Q} ) $ with Prim$H_n C(\Omega^\infty  K(
 A,I)\\hat  ,\Bbb Z )_{\Bbb Q}$. Since $H_n (CC  (A\\hat \, )[1]_{\Bbb Q} ) $ equals  Prim$H_n \Sym (CC  (A\\hat \, )[1]_{\Bbb Q} )$, we need to show that $\chi$ identifies the primitive parts of the homology. The only problem is that 
 the terms of the segment $C(\fgl (A)\\hat \, )_{(a )\,\Bbb Q} \buildrel{(5.4.2)}\over\longleftarrow C(GL(A)^{(m)}\\hat \,   )_{(a )\,\Bbb Q } $ of the chain that defines $\chi$ are not coalgebras. To solve it,
 pick any $a\ge r \ge n$. By 3.5, 3.7 the embedding $C(\fgl_r (A)\\hat\, )\hra  C(\fgl (A)\\hat \,)_{(a)}$ yields isogenies between the homology in degrees $\le n$, so $\tau_{\le n}\chi$ can be  
computed replacing the above segment  by $C(\fgl_r (A)\\hat\, )_{\Bbb Q}\buildrel{(4.1.1)}\over\longleftarrow C(GL_r (A) \\hat \,   )_{\Bbb Q } $. Now all terms of the chain are coalgebras, the maps are compatible with the coproducts (see Remark in 3.5), and we are done.
\qed
\enddemo

\end
\Refs{}
\widestnumber\key{XXXXX}


\ref\key BEK
\by S.~Bloch, H.~Esnault, M.~Kerz
\paper $p$-adic deformations of algebraic cycles classes
\jour AG 1203.2776
\yr 2012
\endref


\ref\key Cr
\by R.~Crew
\paper Crystalline cohomology of singular varieties
\inbook Geometric Aspects of Dwork Theory 
\publ de Gruyter 
\yr 2004
\pages 451--462
\endref




\ref\key DGM
\by B.~Dundas, T.~Goodwillie, R.~McCarthy
\book The local structure of algebraic $K$-theory
\publ Springer
\bookinfo Algebra and applications, vol.~18
\yr 2013
\endref


\ref\key G
\by T.~Goodwillie
\paper Relative algebraic $ K$-theory and cyclic homology
\jour Ann.~Math.
\vol 124
\yr 1986
\pages 347--402
\endref

\ref\key GV
\by A.~Grothendieck, J.-L.~Verdier
\paper 
Prefaisceaux
\inbook 
  Th\'eorie des topos et cohomologie \'etale de sch\'emas (SGA 4), Tome 1
\bookinfo  Lect.~Notes in Math.~269
\publ Springer-Verlag
\yr 1972
\pages 1--218
\endref



\ref\key HKN \by A.~Huber, G.~Kings, N.~Naumann\paper Some complements to the Lazard isomorphism \jour Compositio Math.\vol 147  \yr 2011 \pages 235--262 \endref

\ref\key FI
\by H.~Fausk, D.~Isaksen
\paper t-model structures
\jour Homology Homotopy and Applications
\vol 9 
\yr 2007
\pages 399--438
\endref

\ref\key L
\by M.~Lazard
\paper Groupes analytiques $p$-adiques
\jour Publ.~Math.~IHES
\vol 26
\yr 1965
\pages 5--219
\endref

\ref\key Lo
\by J.-L.~Loday
\book Cyclic homology
\publ Springer
\bookinfo second edition. Grundlehren der mathematischen Wissenschaften
\vol 301
\yr 1998
\endref

\ref\key Lu1 \by J.~Lurie \book Higher topos theory  \bookinfo  Annals of Mathematics Studies \publ Princeton University Press \vol 170 \yr 2009  \endref

\ref\key Lu2 \by J.~Lurie \book Higher algebra  \bookinfo http://www.math.harvard.edu/~lurie/ \yr 2012  \endref

\ref\key M \by J.~P.~May \book The geometry of iterated loop spaces   \bookinfo  Lecture Notes in Math \publ Springer \vol 271 \yr 1972  \endref

\ref\key Mi
\by B.~Mitchell
\paper The cohomological dimension of a directed set
\jour Canad.~J.~Math. \vol 25 \yr1973 \pages 233--238 
\endref

\ref\key NN
\by J.~Nekovar, W.~Niziol  \paper Syntomic cohomology and regulators for varieties over $p$-adic fields \jour  AG 1309.7620 \yr 2013\endref

\ref\key S \by A.~Suslin \paper Stability in algebraic $K$-theory \inbook Algebraic $K$-theory, Oberwolfach 1980  \bookinfo  Lecture Notes in Math \publ Springer \vol 966 \yr 1982 \pages 304--333 \endref

\ref\key TT \by R.~W.~Thomason, T.~Trobaugh \paper Higher algebraic $K$-theory of schemess and of derived categories \inbook Grothendieck Festschrift vol.~III \publ Birkh\"auser \bookinfo Progress in Mathematics \vol 88 \yr 1990 \pages 247--435 \endref



\endRefs
